\documentclass[11pt]{article}
\usepackage{ulem}
\usepackage[makeroom]{cancel}
\usepackage{color}
\usepackage{xcolor}
\usepackage{amsfonts}
\usepackage{amssymb}
\usepackage{amscd}
\usepackage{amsmath}
\usepackage{epsfig}
\usepackage{graphicx, subfigure}
\usepackage{latexsym}
\usepackage{pst-3dplot}
\usepackage{verbatim}
\usepackage{tikz}
\usepackage{pdfsync}
\usepackage{hyperref}
\usepackage{multirow}
\usepackage[utf8]{inputenc}

\usepackage{color}
\usepackage{graphicx}

\usepackage{caption}
\usepackage{lipsum}
\newlength{\dinwidth}
\newlength{\dinmargin}
\newtheorem{theorem}{Theorem}

\newtheorem{proposition}{Proposition}
\newtheorem{corollary}{Corollary}
\newtheorem{remark}{Remark}
\newtheorem{lemma}{Lemma}
\newtheorem{example}{Example}

\def \surf{\mathcal C_\Lambda}
\def\surfx{\mathcal C_\x}
\def\proj{\hat{\mathcal C}}

\def\max{\mathrm{max}}

\def \i{{\rm i}}

\def\pxi{P_{x_i}}
\def\pxj{P_{x_j}}
\def\pxn{P_{x_n}}
\def\pxk{P_{x_k}}
\def\pum{P_{u_m}}
\def\puj{P_{u_j}}
\def\pui{P_{u_i}}

\def\puk{P_{u_k}}

\def\pinfty{{P_\infty}}

\def\S{\mathcal S}
\def\Sx{\mathcal S_\x}

\def\l{\lambda}

\def\plj{P_{\l_j}}

\def\plk{P_{\l_k}}

\def\a{{ a}}
\def\b{{ b}}
\def\x{{\bf x}}

\def\X{{\cal X}}

\def\P{\mathbb P}

\makeatletter
\@addtoreset{equation}{section}
\makeatother

\begin{document}
\title{Deformations of the Hill curves  and isoperiodicity in the KdV and the sine-Gordon equations }

\author{Vladimir Dragovi\'c$^1$ and Vasilisa Shramchenko$^2$}
\date{}

\maketitle

\footnotetext[1]{Department of Mathematical Sciences, The University
	of Texas at Dallas, 800 West Campbell Road, Richardson TX 75080,
	USA. Mathematical Institute SANU, Kneza Mihaila 36, 11000
	Belgrade, Serbia.  E-mail: {\tt
		Vladimir.Dragovic@utdallas.edu}--the corresponding author}

\footnotetext[2]{Department of mathematics, University of
	Sherbrooke, 2500, boul. de l'Universit\'e,  J1K 2R1 Sherbrooke, Quebec, Canada. E-mail: {\tt Vasilisa.Shramchenko@Usherbrooke.ca}}

\

{\it Dedicated to the memory of Sergei Petrovich Novikov (1938-2024).}

\

\begin{abstract}
We consider a family of genus $g$ hyperelliptic curves as double  ramified coverings over the Riemann sphere with the set  of  branch points of the form  $\{0, \infty, x_1, \dots, x_g, u_1, \dots, u_g\}$.  The branch point at infinity $P_\infty$ is selected to be a marked point on the Riemann surfaces.  A meromorphic differential $\Omega$ with a unique pole being of order two at $P_\infty$, is completely defined by the values of half of its periods, the $a$-periods. Fixing values of $a$-periods of $\Omega$, we then find a continuous subfamily in the considered family of  hyperelliptic curves along which all the periods of $\Omega$ are constant. This subfamily is defined by the functions $u_j(x_1, \dots, x_g)$, while $x_1, \dots, x_g$ are independent parameters.  We derive a system of differential equations for the functions $u_j(x_1, \dots, x_g)$, which,  remarkably, has rational coefficients. We call this subfamily the {\it isoperiodic deformations} of the hyperelliptic  curves relative to the given differential of the second kind $\Omega.$  We deduce necessary and sufficient conditions for the existence and uniqueness of isoperiodic deformations. We discuss reality conditions as well. Using the obtained results, we solve the following problem for the Korteweg-de Vries and sine-Gordon equations: {\it starting from an algebro-geometric data which generate a real periodic solution of a period $T$, how to deform the data, so that the associated solutions remain periodic with the same period $T$.}
\end{abstract}

\


\smallskip


\newpage
\tableofcontents

\

\section{Introduction}

The integrals of Abelian differentials over closed contours on Riemann surfaces, called {\it periods} of the differentials, is an important concept in Riemann sufrace theory. Fixing a canonical homology basis of a surface, and fixing periods of a differential corresponding to half of the basis (normalization of the differential), the other half becomes an important tool for analysing the structure of the surface. Periods of holomorphic normalized differentials are combined in the {\it Riemann matrix}, the central piece of the theory. Remarkably, periods of normalized meromorphic differentials are an important component of algebro-geometric solutions to equations of mathematical physics.  In particular, periods of normalized meromorphic differentials of the second kind with a unique pole at a branch point, play a significant role  in solving the Korteweg-de Vries (KdV), sine-Gordon, Kadomtsev-Petviashvili equations, the associated integrable hierarchies, as well as in the expressions for finite gap potentials of the Schr\"odinger operator \cite{5, DMN, Dub1981, ItMa, ItMa1, Nov9174, NovDub, Nov}. Certain values of the periods, or more precisely, the Riemann surfaces they are realized on, correspond to periodic solutions.  In this paper we focus on meromorphic differentials of the second kind on the hyperelliptic curves, with a unique pole being of order two at the branch point at infinity.

 In 1895, Korteweg and de Vries derived cnoidal waves as solutions of the  KdV equation expressed in terms of elliptic functions and periodic in the spatial variable $x$ \cite{Nov}. Eighty years later,  very sophisticated algebro-geometric methods were developed for construction of periodic and quasi-periodic solutions of soliton equations,  see e.g. \cite{Nov, Dub1981} and references therein. These methods are based mostly on the so-called isospectral deformations, Lax representations, and the Baker-Akhiezer functions. Along these lines, tremendous success has been achieved, including the solution of the Novikov conjecture, which resolved the classical Riemann-Schottke problem in algebraic geometry, see \cite{Shi}, \cite{Mu}. Yet, for the KdV equation, as well as for many other soliton equations, the problem of singling out the periodic in $x$ solutions of a given period among the quasi-periodic ones has not been effectively solved till nowadays, see e.g. \cite{Nov}, p. 266.

This motivates our present research, where we approach the problem of constructing periodic solutions in the following way: {\it starting from an algebro-geometric data which generate a periodic solution of a period $T$ to a given soliton equation, how to deform the data, so that the deformed solutions remain periodic with the same period  $
\,T$.} This leads us to the study of {\it isoperiodic deformations}, as an antithesis to the isospectral defformations. Reality  of solutions plays a significant role, in particular from the point of view of applications, for example, in physics. Thus, we pay a special attention to the reality conditions and apply them in  the framework of  the isoperiodic deformations of the KdV and sine-Gordon equations.

  On the other hand, meromorphic differentials of the second kind with real periods were used in \cite{GKCalogeroMoser} to define and study foliations of moduli spaces of Riemann surfaces of a fixed genus with one marked point and a jet of local parameters at the marked point. A leaf of such a foliation is a set of triples $(\mathcal L, P, z)$ where $\mathcal L$ is a genus $g$ Riemann surface, $p$ is a marked point on $\mathcal L$ and $z$ is a local coordinate at $p$ and such that for each triple the periods of a chosen meromorphic differential with a unique second order pole at $p$ are the same. These foliations were related in \cite{GKCalogeroMoser} to Calogero-Moser integrable system.  A combinatorial approach to studying these foliations was proposed in \cite{KLS}.

The question we explore in this paper is closely related to such foliations of the moduli space of surfaces. We consider a family of genus $g$ hyperelliptic coverings ramified over the set  of the form $\{0, \infty, x_1, \dots, x_g, u_1, \dots, u_g\}$ in the complex sphere. This set of {\it branch points} defines the corresponding hyperelliptic covering completely.  All coverings being ramified at the point at infinity $P_\infty$, this point becomes the marked point on our Riemann surfaces.  Our meromorphic differential $\Omega$ with its unique pole being the pole of order two at $P_\infty$, is completely defined by the values of half of its periods, the $a$-periods. As the surface varies, we fix the $a$-periods to be some given constants. Then, the other half, the $b$-periods, are transcendental functions of the branch points $\{x_j, u_j\}_{j=1}^g$. Given a set of $a$-periods, we then find a continuous subfamily in the set of all considered hyperelliptic coverings for which all the periods of $\Omega$ are constant. This subfamily is defined by the functions $u_j(x_1, \dots, x_g)$, that is we have $g$ independent parameters $x_1, \dots, x_g$.  We derive a system of differential equations for the functions $u_j(x_1, \dots, x_g)$.  Remarkably, the obtained system of differential equations has rational coefficients. We call this subfamily the {\it isoperiodic deformations} of the hyperelliptic  curves relative to the given differential $\Omega.$   We deduce necessary and sufficient conditions for the existence and uniqueness of isoperiodic deformations.

As an important particular case of our deformations, we obtain deformations of the Hill curves, the spectral curves of the Schr\"odinger operator with periodic real potential. These deformations preserve the period of the corresponding potential while the spectrum of the Schr\"odinger operator is governed by the obtained differential equations.  For classical results about Hill curves see e.g. \cite{MM, Mc, MM1, MT1, MT2}.

A similar question was studied in \cite{GS}. However,
 our approach, the obtained equations, and derived consequences  are completely different from \cite{GS}. Our main technical   tools for doing analysis on Hurwitz spaces of families of hypereliptic curves  are the fundamental Riemann bidifferential $W$ on Riemann surfaces and the Rauch variational formulae. We present this in the next Section \ref{sec:tech}, see in particular Sections \ref{sect_diff} and \ref{sect_Rauch}.

The paper is organized as follows. In Section \ref{sect_hyper}, we give the necessary background on families of  hyperelliptic curves, meromorphic differentials and their variations over such families; we also give a definition of Hill curves.  Section \ref{sect_deformation} is devoted to definition of isoperiodic deformations for hyperelliptic surfaces relative to a differential of the second kind with a unique pole (of order two) at the ramification point at infinity of the curves. In Section \ref{sect_genusone} we consider isoperiodic deformations in genus one and obtain an ordinary differential equation for one of the two non-constant branch point as a function of another one which describes such deformations; we prove their existence and uniqueness. We then discuss the significance of isoperiodic deformation in the theory of one-gap potentials with applications to planar Neumann system and cnoidal waves.   In Section \ref{sect_genusg}, we derive a system of differential equations for $g$ dependent branch points as functions of $g$ independently varying branch points in order for the corresponding hyperelliptic curves to form an isoperiodic family relative to the chosen differential of the second kind. We prove the existence and unicity of isoperiodic families. Sections \ref{sec:twogap} - \ref{sect_sG} are devoted to discussing applications of our isoperiodic deformations of hyperelliptic curves in the theory of integrable systems. In particular, in Section \ref{sec:twogap}, we construct deformations of finite-gap periodic potentials preserving their period.
 In Section \ref{sec:neum},  we discuss periodic solutions to general Neumann problem and their isoperiodic deformations.  In Section \ref{sect_kdv}, we describe deformations preserving the wavevector of  solutions to the KdV equation constructed from hyperelliptic curves as well as isoperiodic deformations of real periodic finite-gap solutions of the KdV equation. Section \ref{sect_sG} is devoted to continuous families of real periodic solutions to the sine-Gordon equation for which the period is the same, that is constant along the family.  In Section \ref{sec_CR} we relate our isoperiodic deformations to the deformations of comb regions preserving the base of the comb while varying the lengths of the vertical slits.

\

\section{Family of hyperelliptic surfaces}\label{sec:tech}

\subsection{Hyperelliptic curves}
\label{sect_hyper}
Consider a general family of hyperelliptic algebraic curves defined by
\begin{equation}
\label{hyper}
\mathcal C = \{(\l,\mu)\in\mathbb C^2 \;|\; \mu^2=\prod_{j=1}^{2g+1}(\l-\l_j)\},
\end{equation}
parametrized by complex variables $\lambda_1, \dots, \lambda_{2g+1}$ that may vary without coinciding. The compact projective curve $\hat{\mathcal C}\subset \mathbb CP^2$ obtained from a curve of this family by adding a point at infinity $[\lambda:\mu:z]=[0:1:0]$ is singular at the point at infinity.
The  associated compact hyperelliptic Riemann surface, which we denote by $\surf$ with $\Lambda=(\lambda_1, \dots, \l_{2g+1}),$ is of genus $g$ and is such that there is a holomorphic map $\pi:\surf\to \mathbb CP^2$ where $\pi(\surf)=\hat{\mathcal C}$ and
\begin{equation}
\label{pi}
\pi:\surf\setminus\pi^{-1}([0:1:0])\to\proj\setminus\{[0:1:0]\}.
\end{equation}
is a biholomorphism. Using this application, and given that the algebraic curves \eqref{hyper} are non-singular, we identify the points of the Riemann surface with those of the projective curve $\proj$ or of algebraic curve $\mathcal C$ as follows. We denote the point at infinity of $\surf$ by $P_\infty$ where $P_\infty=\pi^{-1}([0:1:0])$ and we say that $P=(\lambda, \mu)$ is a point of the surface $\surf$ if $P=\pi^{-1}([\lambda:\mu:1]).$

We now extend the function $\lambda$, defined naturally on the algebraic curve $\mathcal C \to \mathbb C$ as a projection to the first coordinate, to the Riemann surface $\lambda:\surf\to \mathbb CP^1$ by setting $\lambda(P_\infty)=\infty.$ This allows us to represent the hyperelliptic Riemann surface  as the surface of the two-fold ramified covering $(\surf, \lambda)$ ramified over the points $\lambda_1, \dots, \lambda_{2g+1}, \infty, $ called the branch points of the covering. We denote the corresponding ramification points by $P_{\lambda_j}=(\lambda_j, 0)$ for $j=1, \dots, 2g+1$ and $P_\infty.$

The structure of this ramified covering induces the {\it standard local coordinates} on the surface $\surf$ as follows:
\begin{align}
\label{coordinates}
&\zeta_{\l_k}(P)=\sqrt{\l(P) - \l_k} \quad\mbox{if}\quad P\sim P_{\l_k}, \quad\mbox{with}\quad k=1, \dots, 2g+1\,,
\nonumber
\\
& \zeta_\infty(P)= \frac{1}{\sqrt{\l(P)}}\quad\mbox{if}\quad P\sim P_\infty,
\\
& \zeta_Q(P)=\l(P)-\l(Q) \quad\mbox{if}\quad P\sim Q \quad\mbox{and $Q$ is a regular point.}
\nonumber
\end{align}
It is with these standard local charts that we work in this paper. The variation of $\lambda_j$ induces a change of the curve in the family \eqref{hyper} and of the associated Riemann surface. We thus have a family of compact Riemann surfaces $\surf$ associated with the family of algebraic curves \eqref{hyper}.
%Alternatively, this variation may be considered as a variation of the complex structure of the Riemann surface of the two-fold ramified covering of $\mathbb CP^1$ just defined.

The set of branch points $\Lambda=(\lambda_1, \dots, \lambda_{2g+1})$ varies in  $\mathbb C^{2g+1}\setminus \Delta$ where $\Delta$ is a union of all the diagonals, $\Delta=\{(\lambda_1, \dots, \lambda_{2g+1})\; |\; \lambda_{i}=\lambda_j \text{  for } i\neq j\}.$
We work with a connected domain $\S$ in $\mathbb C^{2g+1}\setminus \Delta$ such that for the hyperelliptic compact Riemann surfaces from the family $\surf$ with $\Lambda\in \S$, we may choose a canonical homology basis $a_1, \dots, a_g; b_1, \dots, b_g$ such that the images $\lambda(a_j) $ and $\lambda(b_j)$ of the basis cycles are fixed for all $\Lambda \in \S$, that is are the same for all surfaces from the family.

\subsection{Abelian differentials}
\label{sect_diff}
We are going to work with Abelian differentials, that is meromorphic differentials, on the compact Riemann surfaces $\surf$ corresponding to  hyperelliptic curves from the family \eqref{hyper}.  In this case, many Abelian differentials may be expressed in terms of the functions $\lambda(P)$ and $\mu(P)$ on the surface, after their extension to the point at infinity by $\lambda(P_\infty)=\infty$ and $\mu(P_\infty)=\infty$ as in Section \ref{sect_hyper}. On such a compact hyperelliptic Riemann surface of genus $g$, we have a basis in the space of holomorphic differentials given by $\phi(P), \;\lambda(P)\phi(P), \dots, \lambda(P)^{g-1}\phi(P)$ with
\begin{equation}
\label{phi}
\phi(P) = \frac{d\lambda(P)}{\mu(P)}.
\end{equation}
Another basis in this space   is given by the holomorphic {\it normalized} differentials $\omega_1, \dots, \omega_g$ satisfying the following conditions with respect to the chosen, as in Section \ref{sect_hyper}, canonical homology basis $a_1, \dots, a_g; b_1, \dots, b_g$ of the surface:
\begin{equation}
\label{normalization}
\oint_{a_j}\omega_k=\delta_{jk}, \qquad j,k=1, \dots, g.
\end{equation}
The {\it $b$-periods} of this basis define the {\it Riemann matrix} $\mathbb B$ of the surface:
\begin{equation}
\label{B}
\oint_{b_j}\omega_k=\mathbb B_{jk}, \qquad j,k=1, \dots, g.
\end{equation}
%which is symmetric and has a strictly positive imaginary part.

\begin{example}
\label{ex:realhyp}
 If a hyperelliptic curve \eqref{hyper} is given by an equation with all $\lambda_j$  real, then it admits an antiholomorphic involution
$$
\rho(\lambda, \mu)=(\bar \lambda, \bar\mu),
$$
see e.g. \cite{DubNat}; the canonical basis of cycles can then be chosen so that  $\rho a_j=a_j$, $\rho b_j=-b_j$, $j=1, \dots, g$. In this case, the holomorphic normalized differentials
$$\omega_j=\frac{\sum_{k=0}^{g-1}c_{jk}\lambda^k}{\mu}, \, j=1,\dots, g$$
have all the coefficients $c_{jk}$ real. The Riemann matrix $\mathbb B$   is purely imaginary in this case,  see e.g. \cite{DubNat}. More generally, a hyperelliptic curve  of equation of the form $y^2=\mathcal P(\lambda)$ where $\mathcal P$ is any real polynomial (having possibly pairs of complex roots) also admits the antiholomorphic involution. These are two classes of {\it real} hyperelliptic curves, where, by definition, real are those curves that admit an antiholomorphic involution. The former class is going to be used several times later, while the latter will be mentioned with respect to real solutions of the sine-Gordon equation.

 {\it A real oval} of a real curve is a connected component of the fixed points of the involution $\rho$. If a curve is given by equation \eqref{hyper} with real  $\lambda_j$ such that $\lambda_1<\lambda_2<\dots<\lambda_{2g+1}$, then it has $g+1$ real ovals $\mathcal O_j=\{(\lambda, \mu)| \lambda_{2j-1}\le \lambda \le \lambda_{2j}, \mu = \pm\sqrt{\prod_{k=1}^{2g+1}(\lambda-\lambda_k)}\}$, $j=1, \dots, g+1$, where $\lambda_{2g+2}=\infty$.

If a real  hyperelliptic curve is defined by a real polynomial, then at least one of $\lambda_j$ is real. The largest real $\lambda_j$ belongs then to the real oval passing through the point at infinity.
\end{example}

An important tool for working with Abelian differentials is the so-called Riemann fundamental bidifferential $W(P,Q)$ depending on the Riemann surface, a choice of a canonical homology basis, and two points $P$ and $Q$ of the surface. It can be expressed in terms of an odd theta-function or defined by its three characteristic properties as follows:
\begin{itemize}
\item Symmetry: $W(P,Q) = W(Q,P);$
\item A second order pole along the diagonal $P=Q$ with the following local expansion in terms of a local parameter $\zeta$ near $P=Q$ and no other singularities:
\begin{equation}
\label{W}
W(P,Q) \underset{P\sim Q}{=} \left( \frac{1}{(\zeta(P) - \zeta(Q))^2}  + {\cal O}(1) \right)d\zeta(P) d\zeta(Q);
\end{equation}
\item Normalization by vanishing of the $a$-periods: $\oint_{a_j} W(P,Q) = 0$ for all $j=1, \dots, g.$ Due to the symmetry,  these integrals can be computed with respect to either $P$ or $Q$.
\end{itemize}
This definition can be shown to imply
\begin{equation}
\label{bW}
\oint_{b_j}W(P,Q) = 2\pi{\rm i}\,\omega_j(P), \qquad j=1, \dots, g.
\end{equation}
We are also going to work with Abelian differentials {\it evaluated} at ramification points of the covering $(\surf, \lambda).$ The evaluation is defined through the standard local parameters from Section \ref{sect_hyper} as follows. Let  $\Upsilon$ be an Abelian differential on $\surf\,.$ We define its {\it value} at a point $\tilde Q\in\surf$ which is not a pole of $\Upsilon$ as the constant term of the Taylor series expansion of the differential with respect to the standard local parameter $\zeta$ from the list \eqref{coordinates} at $\tilde Q$, that is
\begin{equation}
\label{evaluation}
\Upsilon(\tilde Q) = \frac{\Upsilon(P)}{d\zeta(P)}\Big{|}_{P=\tilde Q}.
\end{equation}
As an example, for the value at a ramification point $P_{\lambda_j}$ of the holomorphic differential \eqref{phi}, we use the local parameter $\zeta_{\l_j}(P) = \sqrt{\l(P)-\l_j}$ for which we have $d\lambda(P)=2\zeta_{\l_j}(P)d\zeta_{\l_j}(P)$.
Thus
\begin{equation}
\label{phi-eval}
\phi(\plj) =\frac{\phi(P)}{d\zeta_{\l_j}(P)}\Big{|}_{\zeta_{\l_j}=0}= \frac{2}{\sqrt{\prod_{\substack{i=1\\i\neq j}}^{2g+1}(\l_j-\l_i)}}.
\end{equation}
This definition of evaluation is extended to the bidifferential $W$ by considering one of the arguments of $W$ fixed and thus obtaining a differential with respect to the other argument. For example,
\begin{equation}
\label{W-eval}
W(P, P_\infty) =\frac{W(P, Q)}{d\zeta_\infty(Q)}\Big{|}_{\zeta_{\infty}=0}.
\end{equation}
In this paper, the main role is played by the Abelian differential of the second kind defined by
\begin{equation}
\label{Omega}
\Omega(P)=W(P,P_\infty)+\alpha\omega(P),
\end{equation}
where by $\omega$ we denote the column vector of holomorphic normalized differentials \eqref{normalization}, $\omega=(\omega_1, \dots, \omega_g)^T,$ and $\alpha=(\alpha_1, \dots, \alpha_g)\in \mathbb C^g$ is a constant row vector; its components define the $a$-periods of $\Omega.$ Differential \eqref{Omega} has a pole of the second order at the point $P_\infty$ at infinity of $\surf$, the compact hyperelliptic surface corresponding to a curve from the family \eqref{hyper} and no other singularities.

\subsection{Rauch variation}
\label{sect_Rauch}

Abelian differentials depend on the complex structure of the Riemann surface they are defined on. For our family of surfaces, the complex structure on each of them is given by the local charts \eqref{coordinates}. These complex structures are thus determined by the values of the finite branch points $\lambda_1, \dots, \lambda_{2g+1}$. In this section, we describe a way to study the dependence of the Abelian differentials on  small variation of the $\lambda_j.$

We only consider variations of $\Lambda=(\lambda_1, \dots, \lambda_{2g+1})$ that leave $\Lambda$ in $\S$ defined at the end of Section \ref{hyper}, so that the chosen canonical homology basis does not change under the variation.  An Abelian differential $\Upsilon(P)$ depends on the Riemann surface $\surf$ determined by $\Lambda$ and on a point $P$ on $\surf.$ The Rauch variational formulas \cite{Fay92}  give the dependence of differentials $\Upsilon(P)$ on $\Lambda$ provided the point $P$ is fixed by the condition $\lambda(P)={\rm const}.$ Let us call this derivation the {\it Rauch derivative}:
\begin{equation}
\label{Rauch-der}
\frac{\partial^{{\rm Rauch}}}{\partial \l_k} \Upsilon(P) := \frac{\partial}{\partial \l_k}\Big{|}_{\lambda(P)=const} \Upsilon(P)\,.
\end{equation}
 In the case of the Riemann bidifferential we need to require that both $\l(P)$ and $\l(Q)$ stay fixed. We have the following Rauch variational formula for the  $W$, see \cite{Fay92, KokoKoro}:
\begin{equation*}
\frac{\partial^{{\rm Rauch}}}{\partial \l_k}W(P,Q) := \frac{\partial}{\partial \l_k}\Big{|}_{\substack{\lambda(P)=const\\\lambda(Q)=const}}W(P,Q) = \frac{1}{2}W(P, \plk)W(\plk, Q).
\end{equation*}
This variation of the Riemann bidifferential implies the following Rauch formulas for $\omega_j$ due to \eqref{bW} and for the Riemann matrix  $\mathbb B$ \eqref{B}:
\begin{equation}
\label{RauchB}
\frac{\partial^{{\rm Rauch}} \omega_j(P)}{\partial \l_k} = \frac{1}{2}\omega_j(\plk) W(P, \plk),
\qquad\qquad
\frac{\partial^{{\rm Rauch}} \mathbb B_{ij}}{\partial \l_k} = \pi \i\,\omega_j(\plk)\omega_i(\plk)\,.
\end{equation}
%
%Taking $P=P_s$ for some $s\neq k$ in the first formula in \eqref{RauchB} gives the derivative of $\omega(P_s)$ with respect to $\lambda_k.$
These Rauch variational formulas applied to definition \eqref{Omega} of our main Abelian differential of the second kind $\Omega$, yield
\begin{equation}
\label{Rauch-Omega}
\frac{\partial^{{\rm Rauch}} \Omega(P)}{\partial \l_k} = \frac{1}{2}\Omega(\plk) W(P, \plk)\,.
\end{equation}
Note that in the above formulas we can evaluate both sides of the equation at $P=\plj$ for some $j\neq k,$ obtaining, for example
\begin{equation}
\label{Rauch-Omegalj}
\frac{\partial^{{\rm Rauch}} \Omega(\plj)}{\partial \l_k} = \frac{1}{2}\Omega(\plk) W(\plj, \plk)\,.
\end{equation}
However, this method does not allow us to express the derivative of $\Omega(\plk)$ with respect to $\l_k$ since, for example, that would require having in the right hand side of \eqref{Rauch-Omegalj} the bidifferential $W$ evaluated twice at the same ramification point. This cannot be done as $W$ has a pole when both arguments coincide. In order to circumvent this difficulty, we prove the following result.
\begin{lemma}
\label{lemma_epsilon}
Let $\surf$ be the compact Riemann surfaces associated with the hyperelliptic curves of the family \eqref{hyper}. Let $\plj, \;j=1, \dots, 2g+1,$ be the ramification points of the two-fold ramified covering $\lambda:\surf\to\mathbb CP^1.$ Let $\Omega$ be the meromorphic differential of the second kind on $\surf$ defined by \eqref{Omega} for some fixed $\alpha\in\mathbb C$ and $\Omega(\plj)$ its values at the ramification points defined with respect to the standard local parameters $\zeta_{\l_j}(P)=\sqrt{\lambda(P)-\lambda_j}$ as in \eqref{evaluation}. Then the Rauch derivatives of $\Omega(\plk)$ satisfy
\begin{equation}
\label{epsilon-Omega}
\sum_{\substack{j=1\\j\neq k}}^{2g+1} \frac{\partial^{\rm Rauch}}{\partial \l_j}\Omega(P_{\l_k})+\frac{\partial \Omega(P_{\l_k})}{\partial \l_k}=0.
\end{equation}
\end{lemma}
{\it Proof.} Let $\varepsilon$ be a complex number sufficiently close to zero and let  the ramified covering  $\lambda_\varepsilon: \surf \to \mathbb CP^1$ be defined by $\lambda_\varepsilon: P \mapsto \lambda(P)+\varepsilon\,,$ where $P=(\lambda(P),\mu(P))$ is a point of the surface $\surf$. For $\varepsilon=0$ we obtain the covering $\lambda$ from Section \ref{sect_hyper}.

The covering $\lambda_\varepsilon$  is ramified at the  ramification points $\plj=(\lambda_j,0)$, $\;j=1, \dots, 2g+1,$ which are zeros of $d\lambda_\varepsilon,$ and at $P_\infty=(\infty, \infty)$ for any $\varepsilon$;  the branch points of $\lambda_\varepsilon$ are at $\lambda_1+\varepsilon, \dots, \lambda_{2g+1}+\varepsilon,\, \infty$. %The standard local parameters \eqref{coordinates} are not affected by the shift $\lambda\mapsto \lambda+\varepsilon.$

For the meromorphic differential $\Omega$ defined by \eqref{Omega}  on the surface $\surf$, its value at $\plk$ may be obtained with respect to  the standard local parameters induced  either by the covering $\lambda$ or by the covering $\lambda_\varepsilon$. This could give two quantities that we denote respectively by  $\Omega(\plk)$ and  $\Omega_\varepsilon(\plk)$. However, given that the local parameter $\zeta_{\plk}$ \eqref{coordinates} is unaffected by the shift $\lambda\mapsto \lambda+\varepsilon$, the two quantities coincide: $\Omega_\varepsilon(\plk)=\Omega(\plk).$ This implies in particular $\frac{d}{d\varepsilon}\Omega_\varepsilon(\plk)=0.$

Note that the differential $W(P, P_\infty)$ from definition of $\Omega$ is unaffected by the shift $\lambda\mapsto \lambda+\varepsilon$ as it can be invariantly defined on the Riemann surface $\surf$ as the unique meromorphic normalized ($\oint_{a_n}W(P, P_\infty)=0$ for all $n=1, \dots, g$) differential of the second kind with a pole of second order at $P_\infty$ and  no other singularities.

On the other hand, the quantity $\Omega(\plk)$ is a function of the branch points $\lambda_1, \dots, \lambda_{2g+1}$ of the covering $\lambda:\surf\to\mathbb CP^1$. At the same time,  $\Omega_\varepsilon(\plk)$ is a function of the branch points $\lambda_1+\varepsilon, \dots, \lambda_{2g+1}+\varepsilon$ of the covering $\lambda_\varepsilon$ and thus we have
\begin{equation*}
0=\frac{d}{d\varepsilon}{\Big |}_{\varepsilon=0}\Omega_\varepsilon(\plk)=\sum_{j=1}^{2g+1} \frac{\partial}{\partial (\l_j+\varepsilon)}\Omega_\varepsilon(P_{\l_k}){\Big |}_{\varepsilon=0}=\sum_{j=1}^{2g+1} \frac{\partial}{\partial \l_j}\Omega(P_{\l_k}),
\end{equation*}
which proves the lemma.
$\Box$

This lemma allows us to obtain the derivative of $\Omega(\plk)$ with respect to $\lambda_k$ through the derivatives of $\Omega(\plk)$ with respect to all the $\lambda_j$ with $j\neq k.$ By slight abuse of terminology, we call this derivative the Rauch derivative as well and obtain the following corollary.
\begin{corollary}
\label{cor_epsilon}
In the situation of Lemma \ref{lemma_epsilon}, the following holds
\begin{equation*}
\frac{\partial \Omega(\plk)}{\partial \l_k} = -\frac{1}{2}\sum_{\substack{j=1\\j\neq k}}^{2g+1} W(\plk, \plj)\Omega(\plj)\,.
\end{equation*}
\end{corollary}
{\it Proof.}
This is a direct consequence of Lemma \ref{lemma_epsilon} and the Rauch formulas \eqref{Rauch-Omega}.
$\Box$

Analogously to the proof of Lemma \ref{lemma_epsilon}, we can obtain that the sum of Rauch derivatives of $\omega_s(P_{\l_k})$ with respect to all finite branch points vanishes and thus
\begin{equation}
\label{omega-epsilon}
\frac{\partial \omega_s(\plk)}{\partial \l_k} = -\frac{1}{2}\sum_{\substack{j=1\\j\neq k}}^{2g+1} W(\plk, \plj)\omega_s(\plj),
\end{equation}
for every $k,s=1, \dots, g$.
\subsection { Hill curves}
\label{sect_Hill}
The Hill curves are particular curves in the hyperelliptic family \eqref{hyper}; they are remarkable due to their relationship with the Korteweg-de Vries equation. One can give several equivalent definitions of a Hill curve \cite{MM}, one of them being as follows. A curve of genus $g>0$ from the family \eqref{hyper} is called  Hill curve if $\lambda_j\in\mathbb R$ for all $j=1, \dots, 2g+1$ and there exists a differential of the second kind on the associated compact Riemann surface having a pole of a second order at $P_\infty$ with real leading coefficient $-T$ of the Laurent series at $P_\infty$ and no other singularities and such that its $a$-periods vanish and the $b$-periods are integer multiples of $2\pi\i.$

In our notation, we say that a Hill curve is a curve from \eqref{hyper} with real $\lambda_j$ carrying a differential $\Omega=\Omega_0$ defined by \eqref{Omega} with $\alpha=(0, \dots, 0)$ and such that
\begin{equation}
\label{eq:hill}
\oint_{b_j}\Omega_0=2\pi\i \,\frac{n_j}{T} \quad\text{with}\quad n_j\in\mathbb Z \quad\text{for some}\quad  T>0,  \quad\text{for all}\quad j=1, \dots, g.
\end{equation}
Note that in the notation of \cite{MM}, where ${\rm e}(P)$ is the {\it unit} on the surface \eqref{hyper},  for the differential $\Omega_0$ defined by \eqref{Omega} with $\alpha=0$ and satisfying \eqref{eq:hill}, we have $\Omega_0=-\frac{1}{T}{\rm d\,log\,e}(P).$

\section{Deformations of Hill curves}
\label{sect_deformation}

In this section, we consider isoperiodic deformations of a differential of the second kind over the family of hyperelliptic curves. Deformations of Hill curves is a special case of such deformations.

\subsection{Isoperiodic deformations of a differential of second kind}
\label{sect_iso}

We want to consider continuous variations of the branch points of the curves \eqref{hyper} which leave all periods of the meromorphic differential $\Omega$ \eqref{Omega} invariant. More precisely, the $a$-periods of $\Omega$ are given by the  vector
 $\alpha=(\alpha_1, \dots, \alpha_g)\in \mathbb C^g$ constant by construction. Let us denote by $\beta=(\beta_1, \dots, \beta_g)\in\mathbb C^g$ its $b$-periods:
\begin{equation}
\label{beta}
\oint_{b_k}\Omega(P)=\beta_k.
\end{equation}
Integrating \eqref{Omega}, we have
\begin{equation}
\label{betak}
\beta_k=2\pi\i \omega_k(P_\infty)+\alpha\mathbb B_k,
\end{equation}
with $\mathbb B_k$ being the $k$th column of the Riemann matrix $\mathbb B.$ As we see from the Rauch formulas \eqref{RauchB}, a generic variation of $\lambda_j$ does not guarantee the constancy of $\beta_k$. % give $\frac{\partial \beta_k}{\partial \lambda_j}=0.$
Therefore, in order to  find deformations of the hyperelliptic curve \eqref{hyper} which allow for the periods of the differential \eqref{Omega} to stay constant, we need to reduce the number of degrees of freedom. To this end, we set $\lambda_1=0$ and split the remaining $2g$ finite branch points of our hyperelliptic curves into two sets: $u_1, \dots, u_g$ and $x_1, \dots, x_g$ and allow $u_1, \dots, u_g$ to depend on $x_1, \dots, x_g$. Thus, in the rest of the paper,  our curves are given by
\begin{equation}
\label{hyperx}
\mathcal C_\x = \{(\l,\mu)\in\mathbb C^2 \;|\; \mu^2=\l\prod_{j=1}^g(\l-u_j)\prod_{j=1}^g(\l-x_j)\},
\end{equation}
where only $g$ branch points $x_1, \dots, x_g$ are allowed to vary independently. Here $\x$ stands for $\x=(x_1, \dots, x_g)$ and we denote the compact hyperelliptic Riemann surfaces associated with these curves also by $\mathcal C_\x.$

The domain $\S$ from Section \ref{sect_hyper} is replaced now by a connected domain $\Sx\subset\mathbb C^g\setminus\Delta_\x$ where  $\Delta_\x=\{(x_1, \dots, x_g)\; |\; \exists \,   i\neq j \text{  such that } x_{i}=x_j\}\cup\{(x_1, \dots, x_g)\; |\; x_{i}=u_j \text{  for some } 1\leqslant i, j\leqslant g\}.$ The domain $\Sx$ should be chosen such that for the hyperelliptic compact Riemann surfaces from the family $\surfx$ with $\x\in \Sx$, we may choose a canonical homology basis $a_1, \dots, a_g; b_1, \dots, b_g$ such that the images $\lambda(a_j) $ and $\lambda(b_j)$ of the basis cycles are fixed for all $\x \in \Sx$, that is are the same for all surfaces from the family.

Our objective is to obtain equations for $u_1, \dots, u_g$ as functions of $x_1, \dots, x_g$ that are satisfied when the surface $\surfx$  deforms in a way to keep the quantities $\beta_k$ \eqref{betak} fixed. We refer to these deformations as {\it isoperiodic deformations} of the pair $(\surfx, \Omega)$ with $\Omega=\Omega(\x)$ given by \eqref{Omega} on $\surfx$.
 We also say that the family of pairs $(\surfx, \Omega(\x))$ parametrized by $\x\in\mathcal S_\x$ is an {\it isoperiodic family of hyperelliptic surfaces  relative to differential $\Omega(\x)$}.
In Section \ref{sect_genusone}, we derive such an equation for genus one surfaces while the general case is treated in Section \ref{sect_genusg}.

\section{Isoperiodic deformations in genus one}
\label{sect_genusone}

Here we consider the family of compact Riemann surfaces corresponding to the elliptic curves given by
\begin{equation}
\label{elliptic}
\mathcal C_x = \{(\l,\mu)\in\mathbb C^2 \;|\; \mu^2=\lambda(\lambda-u)(\lambda-x)\},
\end{equation}
with $x$ being an independently varying branch point and $u$ being a function of $x$.  This is the family \eqref{hyperx} in the case $g=1$.
The differential of the second kind $\Omega$ is defined on these surfaces by \eqref{Omega}. The $b$-period of $\Omega$
 is
\begin{equation}
\label{beta}
\beta=2\pi\i \omega(P_\infty)+\alpha\tau,
\end{equation}
 where $\omega$ is the holomorphic normalized differential on the Riemann surface and $\tau$ is its $b$-period. Differentiating this relation
 with respect to $x$ assuming that $u$ to be a function of $x$ using the Rauch formula \eqref{RauchB} with $\tau=\mathbb B_{11}$ and $\omega=\omega_1$, yields
\begin{equation*}
0=\pi\i W(P_\infty,P_x)\omega(P_x)+\pi\i \alpha \omega^2(P_x) + \left(\pi\i W(P_\infty, P_u) +\pi\i \alpha \omega^2(P_u) \right)u',
%\\
%=\pi\i \omega(P_x) \Omega(P_x) + \pi \i \omega(P_u) \Omega(P_u) u',
\end{equation*}
where $u'$ stands for $\frac{du}{dx}$. From this, using the definition \eqref{Omega} of $\Omega$ and the evaluation \eqref{evaluation} at $P_x$ and $P_u$,  we obtain an expression for the first derivative of $u$:
\begin{equation}
\label{uprime}
u'=-\frac{\omega(P_x) \Omega(P_x) }{ \omega(P_u) \Omega(P_u)}.
\end{equation}

\subsection{The equation describing isoperiodic deformations in genus one}
\label{sect_genusone_ode}
Let us now differentiate expression \eqref{uprime} for the first derivative of $u$ the second time with respect to $x$. With the help of Corollary \ref{cor_epsilon} and \eqref{omega-epsilon} for derivatives of differentials evaluated at $P_x$ and $P_u$ with respect to $x$ and $u$, we have
\begin{multline}
\label{u''}
u''=-\frac{\omega(P_x) \Omega(P_x) }{ \omega(P_u) \Omega(P_u)} \left\{ \frac{\frac{d}{dx}\omega(P_x)}{\omega(P_x)} +\frac{\frac{d}{dx}\Omega(P_x)}{\Omega(P_x)}  - \frac{\frac{d}{dx}\omega(P_u)}{\omega(P_u)} - \frac{\frac{d}{dx}\Omega(P_u)}{\Omega(P_u)} \right\}
\\
=-\frac{\omega(P_x) \Omega(P_x) }{2 \omega(P_u) \Omega(P_u)} \left\{   -\frac{\omega(P_0)W(P_x,P_0)}{\omega(P_x)} -\frac{\omega(P_u)W(P_x,P_u)}{\omega(P_x)} +u'\frac{\omega(P_u)W(P_x,P_u)}{\omega(P_x)}\right.
\\
 -\frac{\Omega(P_0)W(P_x,P_0)}{\Omega(P_x)} -\frac{\Omega(P_u)W(P_x,P_u)}{\Omega(P_x)} +u'\frac{\Omega(P_u)W(P_x,P_u)}{\Omega(P_x)}
\\
-\frac{\omega(P_x)W(P_x, P_u)}{\omega(P_u)} + u'\left( \frac{\omega(P_0)W(P_0, P_u)}{\omega(P_u)} +\frac{\omega(P_x)W(P_x, P_u)}{\omega(P_u)} \right)
\\
\left.
-\frac{\Omega(P_x)W(P_x, P_u)}{\Omega(P_u)} + u'\left( \frac{\Omega(P_0)W(P_0, P_u)}{\Omega(P_u)} +\frac{\Omega(P_x)W(P_x, P_u)}{\Omega(P_u)} \right)
\right\}.
\end{multline}
We now replace $W$ evaluated at ramification points by the following expressions in terms of the function $\lambda(P)$ defined on  the Riemann surface $\surfx$ by a natural extension from the algebraic curve \eqref{elliptic} to the point at infinity: $\lambda(P_\infty)=\infty$, see Section \ref{sect_Hill}.
\begin{eqnarray*}
&&W(P, P_0)=\left( \frac{1}{\lambda(P)\omega(P_0)}+I^0\right)\omega(P), \quad W(P, P_x)=\left( \frac{1}{(\lambda(P)-x)\omega(P_x)}+I^x\right)\omega(P),
\\
&&W(P, P_u)=\left( \frac{1}{(\lambda(P)-u)\omega(P_u)}+I^u\right)\omega(P),
\end{eqnarray*}
where $I^0, I^x, I^u$ are normalization constants.
From here, evaluating with respect to the second argument, we obtain
\begin{eqnarray}
\label{Wx0}
&&W(P_x, P_0)=\left( \frac{1}{x\omega(P_0)}+I^0\right)\omega(P_x)=\left( -\frac{1}{x\omega(P_x)}+I^x\right)\omega(P_0),
\\
\label{Wxu}
&&W(P_x, P_u)=\left( \frac{1}{(x-u)\omega(P_u)}+I^u\right)\omega(P_x)=\left( \frac{1}{(u-x)\omega(P_x)}+I^x\right)\omega(P_u).
\\
\label{Wu0}
&&W(P_0, P_u)=\left(- \frac{1}{u\omega(P_u)}+I^u\right)\omega(P_0)=\left( \frac{1}{u\omega(P_0)}+I^0\right)\omega(P_u).
\end{eqnarray}
From \eqref{Wx0} and \eqref{Wxu} we obtain the following relationship between the normalization constants:
\begin{equation}
\label{I-rel}
I^0\!=-\frac{1}{\omega(P_0)x}+\!\left( \! I^x - \frac{1}{x\omega(P_x)}\! \right) \!\frac{\omega(P_0)}{\omega(P_x)} \quad\text{and}\quad
I^u\!=\!\left(\! \frac{1}{(u-x)\omega(P_x)}+I^x \!\right)\!\frac{\omega(P_u)}{\omega(P_x)} - \frac{1}{(x-u)\omega(P_u)}.
\end{equation}
Note that for the holomorphic normalized differential $\omega,$ we can use the following expression in terms of the functions $\lambda$ and $\mu$ on the elliptic curve \eqref{elliptic}:
\begin{equation}
\label{omega-genus1}
\omega(P)=\frac{d\lambda}{I_0\mu} = \frac{d\lambda}{I_0\sqrt{\lambda(\lambda-u)(\lambda-x)}},
\end{equation}
where $I_0$ is the normalization constant depending on the curve.
Evaluating this differential at the ramification points with respect to the standard local parameters according to \eqref{evaluation}, we have
\begin{equation}
\label{omega-evaluation}
\omega(P_0)= \frac{2}{I_0\sqrt{ux}}, \qquad
\omega(P_u)=\frac{2}{I_0\sqrt{u(u-x)}}, \qquad
\omega(P_x)=\frac{2}{I_0\sqrt{x(x-u)}}\,.
\end{equation}
This implies the following identity:
\begin{equation}
\label{omega-identity}
\omega^2(P_0)+\omega^2(P_u)+\omega^2(P_x)=0.
\end{equation}

Now, we  plug \eqref{Wx0}-\eqref{Wu0} as well as \eqref{uprime} in the expression \eqref{u''} for the second derivative of $u$, choosing to use the normalization constant $I^x$ when possible, and then simplify expressions with the help of \eqref{omega-evaluation} and \eqref{omega-identity}. This gives

\begin{multline}
\label{u''-2}
u''
=\frac{u'}{2} \left(  \frac{1}{u} - \frac{1}{x}-\frac{x}{u(x-u)}
+u'\left( \frac{x}{u(x-u)}+\frac{I^x\omega^2(P_u)}{\omega(P_x)}\right) \right)
\\
 + \frac{\omega(P_0) \Omega(P_0) }{2 \omega(P_u) \Omega(P_u)}\left(-\frac{1}{x}+I^x\omega(P_x)\right)  +\frac{1-u'}{2}\left( \frac{1}{u-x}+I^x\omega(P_x)\right)
\\
+\frac{u'}{2}\left(- \frac{1}{u-x} + u'\left(  \frac{1}{u}+I^0\omega(P_0) + \frac{1}{u-x}+I^x\omega(P_x) \right) \right)
+\frac{(u')^2}{2 }\left(\frac{x}{u(x-u)}+\frac{I^x\omega^2(P_u)}{\omega(P_x)}\right)
\\
+\frac{(u')^2}{2 }\frac{\omega(P_0) \Omega(P_0) }{ \omega(P_u) \Omega(P_u)} \left( \frac{x}{u(u-x)}+\frac{I^0\omega^2(P_u)}{\omega(P_0)}\right)
 -\frac{(u')^3}{2 }\left(\frac{x}{u(x-u)}+\frac{I^x\omega^2(P_u)}{\omega(P_x)}\right)\,.
\end{multline}
Let us now consider the following identity on the elliptic curve
\begin{equation*}
\omega(P_0) \Omega(P_0)+\omega(P_u) \Omega(P_u)+\omega(P_x) \Omega(P_x)=0,
\end{equation*}
which can be obtained as the vanishing sum of residues of the differential $\frac{\omega(P) \Omega(P)}{d\lambda(P)}$  on the compact Riemann surface corresponding to the elliptic curve \eqref{elliptic}. Dividing this relation by $\omega(P_u) \Omega(P_u)$ and using  \eqref{uprime} for $u'$, we have
\begin{equation}
\label{Omega-identity}
\frac{\omega(P_0) \Omega(P_0)}{\omega(P_u) \Omega(P_u)}=u'-1.
\end{equation}
Plugging \eqref{Omega-identity} into \eqref{u''-2} and expressing $I^0$ in terms of $I^x$ as in \eqref{I-rel}, we see with the help of \eqref{omega-identity} that the terms with $I^x$ cancel out and thus we obtain an ordinary differential equation \eqref{ode} for $u(x)$ and prove the following theorem.
\begin{theorem}
\label{thm_ode}  Let $\Sx$ be a connected domain as discussed in Section \ref{sect_iso} such that for $x\in\Sx$ one can choose a canonical homology basis $a, b$ on all Riemann surfaces corresponding to the curves from the family
\begin{equation}
\label{subfamily1}
\mathcal C_x = \{(\l,\mu)\in\mathbb C^2 \;|\; \mu^2=\lambda(\lambda-u(x))(\lambda-x)\}
\end{equation}
 in such a way that $\lambda(a)$ and $\lambda(b)$ are fixed for all $x\in\Sx.$
Let $u(x)$ for $x\in\Sx$ be such that the compact Riemann surfaces associated with curves of subfamily \eqref{subfamily1} of the family \eqref{elliptic}
are such that the differential of the second kind $\Omega$ defined on each of these surfaces by \eqref{Omega} has $a$- and $b$-periods that are constant for all $x\in\S$.
%along the subfamily \eqref{subfamily1}.
Then $u(x)$ satisfies the following equation:
\begin{equation}
\label{ode}
u''
=\frac{1}{2}\left( \frac{1}{x}+\frac{1}{u-x}\right)-\frac{u'}{2} \left(   \frac{2}{x}+\frac{1}{u-x} \right)
+\frac{(u')^2}{2 } \left(\frac{2}{u}+\frac{1}{x-u}  \right)
 -\frac{(u')^3}{2 }\left(\frac{1}{u}+\frac{1}{x-u}\right)\,.
\end{equation}
In particular if  the curves \eqref{subfamily1} are Hill curves for every $x\in\S$, then $u(x)$ is a solution of \eqref{ode}.
\end{theorem}
Conversely, some solutions of equation \eqref{ode} correspond to  isoperiodic deformations as explained in the next theorem.
\begin{theorem}
\label{thm_converse-g1}
Let $x_0\in\mathbb C\setminus\{0\}$, and let $\Omega=\Omega(x_0)$ be the meromorphic differential of the second kind defined by \eqref{Omega} on a compact Riemann surface $\mathcal C_{x_0}$ of the elliptic curve \eqref{elliptic} with $x=x_0$.  Assume that $x_0$ is such that $\Omega(x_0; P_u)\neq 0$. Let $\{a,b\}$ be a canonical homology basis on $\mathcal C_{x_0}$ such that the projections $u(a)$ and $u(b)$ do not intersect a certain neighbourhood $\hat \X$ of $x_0$. Then there exists a unique continuous family  $(\mathcal C_x, \Omega(x))$ providing an isoperiodic deformation of the pair $(\mathcal C_{x_0}, \Omega(x_0))$ for $x$ varying in some neighbourhood $\X\subset\hat\X$ of $x_0.$

\end{theorem}
{\it Proof.}
Isoperiodic deformations of $(\mathcal C_{x_0}, \Omega(x_0))$ are defined by the condition $\oint_b\Omega(x)={\rm const},$ as the $a$-period of $\Omega(x)$ is constant by definition \eqref{Omega}.  By the implicit function theorem, the relation $\oint_b\Omega(x)={\rm const}$ defines a function $u(x)$ in a neighbourhood of $x=x_0$ if the derivative $\frac{\partial}{\partial u}\oint_b\Omega(x)$ is non-zero at $x=x_0$. The $b$-period of $\Omega$ is given by \eqref{beta}. Thus, the derivative in question can be calculated by Rauch variational formulas \eqref{RauchB}:
\begin{equation*}
\frac{\partial}{\partial u}\oint_b\Omega = \frac{\partial^{\rm Rauch}}{\partial u}(  2\pi\i \omega(P_\infty)+\alpha\tau  ) = \pi\i W(P_u, P_\infty)\omega(P_u) +\pi\i\alpha\omega^2(P_u) = \pi\i\Omega(P_u)\omega(P_u).
\end{equation*}

The factor of  $\omega(P_u)$ is non-zero since a holomorphic differential on a genus one surface does not vanish, and $\Omega(P_u)$ is non-zero for the curve $\mathcal C_{x_0}$ given the choice of $x_0$. This proves the existence of an isoperiodic deformation in some neighbourhood $\X\subset\hat \X$ of $x_0.$

On the other hand, Theorem \ref{thm_ode} states that for every continuous isoperiodic deformation, the function $u(x)$ satisfies equation \eqref{ode}. In addition,  for an isoperiodic deformation, $u'(x_0)$ is given by \eqref{uprime}. Thus the initial values $u(x_0)$ and $u'(x_0)$ determine a unique solution to \eqref{ode}. This proves the unicity of an isoperiodic deformation of the pair $(\mathcal C_{x_0}, \Omega(x_0)).$ We assume $\X$ small enough so that $u(x)$  does not belong to $u(a)$ or $u(b)$ for $x\in\X.$
$\Box$

\subsection{Isoperiodic deformations of one-gap potentials}
\label{sec:onegap}

 In this section, we discuss an application of Theorems \ref{thm_ode} and \ref{thm_converse-g1} to the theory of one-gap potentials.
Let
$\wp(z)$ be the standard Weierstrass function associated with the elliptic curve
\begin{equation}
\label{GW}
\Gamma_{e_2,e_3}=\{ (\lambda, w)\in\mathbb C^2\,|\,w^2=4\lambda^3-g_2\lambda-g_3\}.
\end{equation}
Denoting the ramification points of the covering $\lambda: \Gamma_{e_2,e_3}\to \mathbb CP^1$ by $P_{e_j}=(\lambda=e_j, w=0)$ and $P_\infty$, the standard choice of the canonical homology basis on $\Gamma_{e_2,e_3}$ is as follows: the cycle $a$ encircles the points $P_\infty$ and $P_{e_1}$ while the cycle $b$ encircles the points $P_\infty$ and $P_{e_2}.$
The elliptic curve can be parameterized by the points of the torus $\mathbb T^2=\mathbb C/\{2m\omega+2n\omega'\}$ with
$\lambda=\wp(z)$ and $w=\wp'(z)$ for $z\in\mathbb T^2$. The Weierstrass functions $\zeta$ and $\sigma$  are defined in the standard manner (see e.g. \cite{Akh1})
by $\zeta'(z)=-\wp(z)$ and $\sigma'(z)/\sigma(z)=\zeta(z)$.

In  1940, Ince \cite{In} showed that $v_1(X)=2\wp(X), \;X\in\mathbb C$ is a one-gap potential for the Lam\'e operator
$$\mathbb L=-\frac{\partial^2}{\partial X^2}+v_1(X).$$
%$$\mathbb L=-\frac{\partial^2}{\partial X^2}+2\wp(X).$$

In \cite{Hoc}, it was shown that there are no other one-gap potentials apart from those constructed by Ince.

Following e.g. \cite{Dub1981} and \cite{Akh1}, let us introduce the Baker-Akhiezer function
$$\psi(X, P)=\frac{\sigma(z-z_1-X)\exp(X\zeta(z))}{\sigma(z-z_1)\sigma(z_1+X)},$$
where $P=(\wp(z), \wp'(z))$. Then, $\psi(X, P)$ is an eigenfunction of the Lam\'e operator $\mathbb L=-\frac{\partial^2}{\partial X^2}+2\wp(X)$:
$$
\mathbb L\psi(X,P)=\lambda \psi(X,P), \quad \lambda = \wp(z).
$$

The Lam\'e potential $v_1(X)$ is a doubly periodic function with periods $2w_1$ and $2w_2,$ where
$$
2w_1=\oint_{a}\phi,\qquad 2w_2=\oint_{b}\phi.
$$
Here $\phi=\frac{d\lambda}{w}$ is a non-normalized holomorphic differential on $\Gamma,$  and $a$ and  $b$ form a canonical basis of homology of $\Gamma$, where we denote the compact Riemann surface associated with the elliptic curve by $\Gamma$ as well.

%Let us pose the following natural question: find deformations of the elliptic curve $\Gamma$ which preserve $2w_1,$ one of the periods of the potential $v_1(X)=2\wp(X)$.

It could thus be interesting to study deformations of the elliptic curve $\Gamma$ which preserve one of the periods of the potential $v_1(X)=2\wp(X)$.

In the standard way, let us denote the zeros of the polynomial $4z^3-g_2z-g_3$ by $e_1,\; e_2,\; e_3$. Then $e_1=-e_2-e_3$.
\begin{theorem}
\label{thm_2w1}
There exists a unique continuous family of elliptic curves \eqref{GW}
parameterized by interdependent values of $e_2,\; e_3 \in \mathbb C$ such that the period $2w_1$ is constant over all curves of the family. Moreover, for this family $u=2e_2+e_3$ satisfies equation \eqref{ode} as a function of $x=e_2+2e_3$.
\end{theorem}
{\it Proof.} Let us perform a shift by $-e_1$ in the $\lambda$-sphere which maps the curves of the family $\Gamma_{e_2,e_3}$ into the curves $\hat\Gamma_{x,u}$ with branch points at $0,\; x=2e_3+e_2,\;u=2e_2+e_3$ and the point at infinity. The canonical homology basis $a$ and $b$ chosen above is shifted accordingly.  The periods  $2w_1, \;2w_2$ of the holomorphic differential $\phi=\frac{d\lambda}{w}$ are unchanged under the change of variables.

Let $\Omega_0$ be the differential of the second kind \eqref{Omega} with $\alpha =0$ defined on the curves from the family $\hat\Gamma_{x,u}$.
Applying \eqref{beta} for $\alpha =0$ we see that the $b$-period of $\Omega_0$ is given by
\begin{equation}
\label{Omega0}
\oint_b\Omega_0=2\pi\i\omega(P_\infty).
\end{equation}
Here $\omega$ is the holomorphic normalized differential on $\hat\Gamma_{x,u}$, that is $\omega=\frac{\phi}{2w_1},$ and the evaluation at the point at infinity of $\hat\Gamma_{x,u}$ is done as in \eqref{evaluation} with respect to the standard local parameter $\zeta(P) = 1/\sqrt{\lambda(P)}:$
\begin{equation*}
\omega(P_\infty) =\frac{\omega(P)}{d\zeta(P)}\Big{|}_{\zeta=0}= \frac{d\left( \frac{1}{\zeta^2}\right)}{2w_1\sqrt{\frac{4}{\zeta^6}-\frac{g_2}{\zeta^2}-g_3}}\Big{|}_{\zeta=0}=-\frac{1}{2w_1}.
\end{equation*}
Thus we obtain that  $2w_1$ is constant for a family $\hat\Gamma_{x,u}$  if and only if the $b$-period of $\Omega_0$ is constant for the same family. Thus the required family of curves is an isoperiodic family $(\hat\Gamma_{x,u}, \Omega_0)$; its existence  and uniqueness is given by Theorem \ref{thm_converse-g1}. The fact that $u(x)$ satisfies equation \eqref{ode}  follows from Theorem \ref{thm_ode}.
$\Box$
\begin{corollary}
The function $u(x)$ from Theorem \ref{thm_2w1} satisfies equation \eqref{ode} with the initial condition given by \eqref{uprime} with $\Omega=\Omega_0$ given by \eqref{Omega} with $\alpha=0$ and $\omega=\frac{d\lambda}{2w_1w}$ being the holomorphic normalized differential on the curve $\hat\Gamma_{x,u}$ for some initial value of $x$.
\end{corollary}
\begin{corollary}
\label{prop:onegap}
There exists a  unique continuous deformation of the curve $\Gamma_{e_2,3_3}$ \eqref{GW} which preserves the period $2w_1$ of the potential $v_1(X)=2\wp(X)$. This deformation is such that $u=2e_2+e_3$ satisfies equation \eqref{ode} as a function of $x=2e_3+e_2$  with the initial condition  expressed by \eqref{uprime} with respect to the following differentials on the curve $\hat\Gamma_{x,u}$ obtained from  $\Gamma_{e_2,e_3}$ by a shift by $-e_1$ in the $\lambda$-sphere: $\Omega=\Omega_0$ given by \eqref{Omega} with $\alpha=0$ and $\omega=\frac{d\lambda}{2w_1w}.$
\end{corollary}
 We now give two examples providing an application of the theory of this section to isoperiodic deformations of solutions to the  planar version of the Neumann system  and of the cnoidal waves expressed as periodic  solutions to the Korteweg-de Vries equation.
\begin{example}
\label{ex:planarneum}
[Period-preserving deformations of solutions of the planar Neumann system]
Let us consider a system of two uncoupled oscillators with frequencies $\sqrt{A_1}$, $\sqrt{A_2}$ with $\sqrt{A_1}< \sqrt{A_2}$, under the assumption that the external force $F$ constrains the particles to remain on a unit circle. The equations of motion are
\begin{align}\label{eq:2neum}
\ddot q_1+A_1q_1&=-Fq_1,\\
\ddot q_2+A_2q_2&=-Fq_2,
\end{align}
where $F=\dot q_1^2+\dot q_2^2-A_1q_1^2-A_2q_2^2$ and $q_1^2+q_2^2=1$ with $q_1,\;q_2$ being the positions of the two particles. We can use the Jacobi elliptic coordinates which map a point $(q_1, q_2)$ from the unit circle to the conical Jacobi coordinate $\mu_1$, such that
$$
\frac{q_1^2}{A_1-\mu}+\frac{q_2^2}{A_2-\mu}=0.
$$
Thus, $\mu_1=A_2q_1+A_1q_2$ satisfies $A_1\le \mu\le A_2$ and
$$
q_1=\sqrt{\frac{\mu-A_1}{A_2-A_1}}, \quad q_2=\sqrt{\frac{A_2-\mu}{A_2-A_1}}.
$$
After differentiation with respect to time, we get a Weierstrass-type differential equation for $\mu$:
$$
\big(\dot \mu\big)^2=-4(\mu-A_1)(\mu-A_2)(\mu-A_1-A_2-2h),
$$
where
$$h=\frac{1}{2}\big(\dot q_1^2+\dot q_2^2+A_1q_1^2+A_2q_2^2)$$
is the first integral of the system equal to the energy of the uncoupled pair of oscillators.
Thus,
$$
\mu(t)=\wp(\i t +t_0)-z_0,\quad z_0=-\frac{2}{3}(h + A_1 + A_2),
$$
with some constant $t_0$, which corresponds to the initial conditions for the planar Neumann system, where $\wp$ is the Weierstrass function of the elliptic curve $\Gamma$ \eqref{GW} with
%$$
 %\Gamma= \{(\l,w)\in\mathbb C^2 \;|\; w^2=4(\lambda-e_1)(\lambda-e_2)(\lambda-e_3)\}, \qquad e_1+e_2+e_3=0.
%$$
 $\lambda=\mu+z_0$ and the standard choice of the canonical homology basis as above.  According to Corollary \ref{prop:onegap},  the  solution $\mu(t)$ is deformed by varying the underlying elliptic curve, in a way that the period $2w_1/\i$ of $\mu(t)$ is preserved if and only if $u=2e_2+e_3$ satisfies  equation \eqref{ode} as a function of $x=e_2+2e_3$ with the initial condition \eqref{uprime},  where $\Omega=\Omega_0$ from the proof of Theorem \ref{thm_2w1}, see \eqref{Omega0}.

The three-dimensional version of the Neumann problem was studied by C. Neumann in 1859. Higher-dimensional generalizations of the Neumann system were studied by Moser \cite{Mo}, \cite{Mo1} and many others. We discuss isoperiodic deformations of solutions to the Neumann system in dimensions three and higher in Section \ref{sec:neum}.
\end{example}

\begin{example}\label{ex:cnoidalwave}[Period-preserving deformations of cnoidal waves]
The function
$$
v(X, t)\equiv v(X-ct)=2\wp(X-ct-X_0)-\frac{c}{6}
$$
provides a periodic  solution to the Korteweg-de Vries equation
$$
v_t-6vv_X+v_{XXX}=0,
$$
 of the form of cnoidal wave. According to Corollary \ref{prop:onegap}, the solution $v(X, t)$ can be continuously deformed by varying the underlying elliptic curve in a way that the period $2w_1$ is preserved. This deformation is such that $u=2e_2+e_3$ satisfies equation \eqref{ode} as a function of $x=e_2+2e_3$   with the initial condition \eqref{uprime},  where $\Omega=\Omega_0$ from \eqref{Omega0}.
\end{example}

\section{Isoperiodic deformations in genus greater than one}
\label{sect_genusg}

\subsection{Variation of branch points}
Consider now the family of curves \eqref{hyperx} for any $g>0$ assuming that the branch points $u_1, \dots, u_g$ are functions of the independently varying branch points $x_1, \dots, x_g.$
With this assumption in mind, we differentiate equation \eqref{betak} with respect to $x_j$ and obtain, in Theorem \ref{thm_main} equations for the functions $\{u_k\}$ with respect to $\{x_j\}.$

Differentiating \eqref{betak} with respect to an independently varying branch point $x_j$ with $j=1, \dots, g$, using Rauch formulas \eqref{RauchB}, we have
\begin{equation*}
0=\frac{\partial \beta_k}{\partial x_j} = \pi\i\omega_{k}(\pxj)(W(\pxj, \pinfty)+\alpha\omega(\pxj)) +\pi \i \sum_{m=1}^g\omega_k(\pum)(W(\pum, \pinfty)+\alpha\omega(\pum)) \frac{\partial u_m}{\partial x_j},
\end{equation*}
which can be rewritten in terms of differential $\Omega$ \eqref{Omega} in the form
\begin{equation}
\label{linsys}
\sum_{m=1}^g\omega_k(\pum)\Omega(\pum) \frac{\partial u_m}{\partial x_j}=-\omega_{k}(\pxj)\Omega(\pxj).
\end{equation}
The latter equation can be seen as a system of $g$ linear equations for $g$ variables $ \frac{\partial u_m}{\partial x_j}$ for every fixed $j.$ In the equations of this system, the differential $\omega_k$ can be replaced by any combination of holomorphic normalized differentials $\omega_1, \dots, \omega_g.$
To obtain a solution to this linear system, it is convenient to introduce a new basis $v_1, \dots, v_g$ in the space of holomorphic differentials defined by the condition
\begin{equation}
\label{vij}
v_j(P_{u_i})=\delta_{ij} \quad\text{ for } i,j=1, \dots, g,
\end{equation}
see \eqref{v} below for a more explicit definition.

Differentials $v_j$ are linear combinations of $\omega_1, \dots, \omega_g$ and thus we can replace $\omega_k$ by $v_k$ in \eqref{linsys}.
Due to relation $v_k(\pum)=\delta_{km}$ this gives the following theorem generalizing \eqref{uprime}.
\begin{theorem}
\label{thm_derivatives}
Let $\Sx$ be such as discussed in Section \ref{sect_iso}. Let $\Omega$ be the differential defined by \eqref{Omega} on the curves of the family $\surfx$ \eqref{hyperx} with $u_j=u_j(x_1, \dots, x_g)$ and the basis of $a$- and $b$-cycles chosen consistently for all curves $\surfx$ with $\x=(x_1, \dots, x_g)\in\Sx$. If the periods of $\Omega$ stay constant for all the curves of the family for $\x\in\Sx$, then derivatives of $u_1, \dots, u_g$ with respect to $x_1, \dots, x_g$ are given by
\begin{equation}
\label{derivative}
 \frac{\partial u_k}{\partial x_j}=-\frac{ v_{k}(\pxj)\Omega(\pxj)}{ \Omega(\puk)}\,,
\end{equation}
where $\pxj$ and $\puk$ are the ramification points corresponding to the branch points $x_j$ and $u_k,$ respectively, and $v_j$ are holomorphic differentials defined by conditions \eqref{vij} on the curves from the family $\surfx$.
\end{theorem}
%
%thus obtaining an expression for derivatives of dependent branch points with respect to the independent ones which generalizes \eqref{uprime}.

Let us note the following useful relation. The sum of residues of the meromorphic differential $\frac{v_m(P)\Omega(P)}{d\l(P)}$ on a curve from the family $\mathcal C_\x$ \eqref{hyper} vanishes and therefore we have
\begin{equation}
\label{Oum}
\Omega(P_0)v_m(P_0) + \sum_{i=1}^g\Omega(\pxi)v_m(\pxi)+\Omega(\pum)=0,
\end{equation}
which, dividing out $\Omega(\pum)$ and using \eqref{derivative} becomes
\begin{equation}
\label{useful}
\frac{\Omega(P_0)v_m(P_0)}{\Omega(\pum)} =  \sum_{i=1}^g\frac{\partial u_m}{\partial x_i}-1.
\end{equation}

\subsection{Coordinate representation of Abelian differentials}

The holomorphic differentials  $v_j$ determined by conditions \eqref{vij} admit the following expressions
in terms of the functions $\lambda$ and $\mu$ on a surface of the family $\mathcal C_\x$ (\eqref{hyper}):
\begin{equation}
\label{v}
v_j(P)=\frac{\phi(P) \prod_{i\neq j} (\l-u_i)}{\phi(P_{u_j}) \prod_{i\neq j} (u_j-u_i)}, \quad j=1, \dots, g,
\end{equation}
%
%%
%\begin{equation}
%\label{v}
%v_j(P)=\frac{\phi(P) \prod_{\substack{i=1\\i\neq j}}^g (\l-u_i)}{\phi(P_{u_j}) \prod_{\substack{i=1\\i\neq j}}^g (u_j-u_i)}, \quad j=1, \dots, g,
%\end{equation}
%%

where $\phi$ is the holomorphic non-normalized differential \eqref{phi} and the index range is $i=1, \dots, g$ in the products.

Analogously to formulas \eqref{Wx0} - \eqref{Wu0} in genus one, we can write expressions for $W$ with one of the two arguments evaluated at a branch point in terms of the function $\lambda(P)$ on the compact surface $\mathcal C_\x$ at the expense of introducing normalization constants denoted by $I_i^{x_k}$ and $I_i^{u_m}$  below.
\begin{eqnarray}
\label{Wppxk}
&&W(P, \pxk) =\frac{\phi(P)}{\phi(\pxk)(\lambda(P)-x_k)} + \sum_{i=1}^g I_i^{x_k}v_i(P),\quad k=1, \dots, g,
\\
\label{Wppum}
&&W(P, \pum) =\frac{\phi(P)}{\phi(\pum)(\lambda(P)-u_m)} + \sum_{i=1}^g I_i^{u_m}v_i(P), \quad m=1, \dots, g.
\end{eqnarray}
Evaluating further, as in \eqref{evaluation}, the second argument at another ramification point, and using definition \eqref{vij} or \eqref{v}  of the differentials $v_j$, we obtain
\begin{eqnarray}
\label{Wpumpxk}
&&W(\pum, \pxk) =\frac{\phi(\pum)}{\phi(\pxk)(u_m-x_k)} +  I_m^{x_k} = \frac{\phi(\pxk)}{\phi(\pum)(x_k-u_m)} + \sum_{i=1}^g I_i^{u_m}v_i(\pxk),
\\
\label{Wpxnpxk}
&&W(\pxn, \pxk) =\frac{\phi(\pxn)}{\phi(\pxk)(x_n-x_k)} + \sum_{i=1}^g I_i^{x_k}v_i(\pxn),
\\
\label{Wpujpum}
&&W(\puj, \pum) =\frac{\phi(\puj)}{\phi(\pum)(u_j-u_m)} + \ I_j^{u_m}.
\end{eqnarray}

Let us now prove five useful technical lemmas.
\begin{lemma} Let $W$ be the Riemann fundamental bidifferential \eqref{W} on a surface $\surfx$ from the family \eqref{hyperx}, let $\phi$ and $v_j$ be the holomorphic differentials \eqref{phi} and \eqref{v}, respectively. Let the differentials be evaluated at ramification points according to \eqref{evaluation} with respect to the standard local parameters \eqref{coordinates}. For any two integers $1\leqslant k,n\leqslant g$,  the following identity holds:
\label{lemma_T1}
\begin{equation}
\label{T1}
\sum_{j=1}^g W(\puj, \pxk)v_j(\pxn) = W(\pxn, \pxk) + \frac{\phi(\pxn)}{\phi(\pxk)} \frac{1}{x_k-x_n} \frac{\prod_{i=1}^g(x_n-u_i)}{\prod_{i=1}^g(x_k-u_i)}\,,
\end{equation}
where $\pxj$ and $\puj$ are the ramification points of $\lambda:\surfx\to\mathbb CP^1$ corresponding to the branch points $x_j$ and $u_j,$ respectively.
\end{lemma}
{\it Proof.} For some fixed $k$ and $n$, consider the following differential on the compact surface $\mathcal C_\x$ having simple poles at $\pxn$, $\pxk$ and $\pui$ with $i=1, \dots, g,$  and no other singularities:
\begin{equation*}
\Upsilon_1(P)=\frac{W(P, \pxk) d\lambda(P)}{\phi(P)(\lambda(P)-x_n) \prod_{i=1}^g(\lambda(P)-u_i)}.
\end{equation*}
The vanishing of the sum of all  residues of this differential  is equivalent to the statement of the lemma. The residues at $\pxn$ and $\pui$ are
\begin{equation*}
\underset{P=\pxn}{\rm res}\Upsilon_1(P)=\frac{2W(\pxn, \pxk) }{\phi(\pxn)\prod_{i=1}^g(x_n-u_i)},
\qquad
\underset{P=\puj}{\rm res}\Upsilon_1(P)=\frac{2W(\puj, \pxk) }{\phi(\puj)(u_j-x_n)\prod_{\substack{i=1\\ i\neq j}}^g(u_j-u_i)}.
\end{equation*}
 Let us compute the residue at $\pxk$, making use of expression \eqref{Wppxk} for $W(P, \pxk)$ and the standard local parameter $\zeta_{x_k}(P)=\sqrt{\lambda(P)-x_k}$ at $\pxk$:
\begin{equation*}
\underset{P=\pxk}{\rm res}\Upsilon_1(P)=\underset{P=\pxk}{\rm res}\frac{ 2\zeta_{x_k}(P) d\zeta_{x_k}(P)\left( \frac{\phi(P)}{\phi(\pxk)\zeta_{x_k}^2(P)} + \sum_{i=1}^g I_i^{x_k}v_i(P)  \right)}{\phi(P)(\lambda(P)-x_n) \prod_{i=1}^g(\lambda-u_i)}
=\frac{ 2}{\phi(\pxk)(x_k-x_n) \prod_{i=1}^g(x_k-u_i)}.
\end{equation*}
The statement of the lemma is now obtained by equating to zero the sum of residues using expression \eqref{v} for $v_j.$
$\Box$
\begin{lemma}
\label{lemma_T2}
Assuming notation of Lemma \ref{lemma_T1}, let $m$ be fixed, $1\leqslant m\leqslant g.$ Let $I_m^{u_m}$  be the normalization constant from \eqref{Wppum}. Then the following identity holds
\begin{equation}
\label{T2}
\sum_{\substack{j=1,\\ j\neq m}}^g W(\puj, \pum)v_j(\pxn) = W(\pxn, \pum)  - \frac{v_m(\pxn)}{x_n-u_m}+ v_m(\pxn)\sum_{\substack{i=1,\\ i\neq m}}^g \frac{1}{u_m-u_i} -v_m(\pxn) I_m^{u_m}\,.
\end{equation}
\end{lemma}
{\it Proof.} The following differential on the compact surface $\mathcal C_\x$ has simple poles at $\pxn$, $\pui$ for $i\neq m,$ a pole of order three at $\pum$ and no other singularities:
\begin{equation*}
\Upsilon_2(P)=\frac{W(P, \pum) d\lambda(P)}{\phi(P)(\lambda(P)-x_n) \prod_{i=1}^g(\lambda(P)-u_i)}.
\end{equation*}
The statement of the lemma is obtained by equating to zero the sum of the residues of $\Upsilon_2.$ Let us compute the residue at $\pum.$ Using \eqref{Wppum} for $W(P, \pum)$ and the local parameter $\zeta_{u_m}(P)=\sqrt{\lambda(P)-u_m},$ we have
\begin{multline*}
\underset{P=\pum}{\rm res}\Upsilon_2(P)
=
\underset{P=\pum}{\rm res} \frac{2d\zeta_{u_m}(P)\left( \frac{\phi(P)}{\phi(\pum)\zeta_{u_m}^2(P)} + \sum_{i=1}^g I_i^{u_m}v_i(P) \right)}{\zeta_{u_m}(P)\phi(P)(\lambda(P)-x_n) \prod_{\substack{i=1\\ i\neq m}}^g(\lambda(P)-u_i)}
\\
=-\frac{2}{\phi(\pum)(u_m-x_n)\prod_{\substack{i=1\\ i\neq m}}^g(u_m-u_i)} \left( \frac{1}{u_m-x_n} + \sum_{\substack{i=1\\ i\neq m}}^g\frac{1}{u_m-u_i} -I_m^{u_m} \right).
\end{multline*}
It remains to compute residues at the simple poles and write the sum of residues using expression \eqref{v} for $v_j(\pxn).$
$\Box$
\begin{lemma}
\label{lemma_T3}
Assuming notation of Lemma \ref{lemma_T1},  let $I_j^{x_k}$ for $j=1, \dots, g$  be the normalization constants from \eqref{Wppxk}. For any integer $k,$  $1\leqslant k\leqslant g,$ the following identity holds:
\begin{equation}
\label{T3}
\sum_{j=1}^g W(\puj, \pxk)v_j(\pxk) =  \sum_{j=1}^g  I_j^{x_k}v_j(\pxk) - \sum_{j=1}^g\frac{1}{x_k-u_j}\,.
\end{equation}
\end{lemma}
{\it Proof.}  This lemma is proved looking at the sum of residues of the following differential on the compact surface $\mathcal C_\x:$
\begin{equation*}
\Upsilon_3(P)=\frac{W(P, \pxk) d\lambda(P)}{\phi(P)(\lambda(P)-x_k) \prod_{i=1}^g(\lambda(P)-u_i)}.
\end{equation*}
This differential has simple poles at $\pui$ for $i=1, \dots, g,$ a pole of order three at $\pxk$ and no other singularities. Its residue at $\pxk$ is obtained writing expression \eqref{Wppxk} for $W(P, \pxk)$ and using $\zeta_{x_k}=\sqrt{\lambda-x_k}$ as a local parameter:
\begin{equation*}
\underset{P=\pxk}{\rm res}\Upsilon_3(P)
%\underset{P=\pxk}{\rm res} \frac{2d\zeta_{x_k}(P)\left( \frac{\phi(P)}{\phi(\pxk)\zeta_{x_k}^2(P)} + \sum_{j=1}^g I_j^{x_k}v_j(P) \right)}{\zeta_{x_k}(P)\phi(P) \prod_{i=1}^g(\lambda(P)-u_i)}
%\\
=-\frac{2}{\phi(\pxk) \prod_{i=1}^g(x_k-u_i)} \sum_{i=1}^g \frac{1}{x_k-u_i}
+
\frac{2\sum_{j=1}^g I_j^{x_k}v_j(\pxk) }{\phi(\pxk) \prod_{i=1}^g(x_k-u_i)}\,.
\end{equation*}
Using again expression \eqref{v} for $v_j,$ we see that the vanishing of the sum of residues of $\Upsilon_3$ is equivalent to the equality claimed in the lemma.
$\Box$

\begin{lemma}
\label{lemma_s}
Let $1\leqslant m\leqslant g.$ Let $I_m^{u_m}$ be the normalization constant from \eqref{Wppum}.
Assuming notation of Lemma \ref{lemma_T1}, with $P_0$ being the ramification point of $\lambda:\surfx\to\mathbb CP^1$ such that $\lambda(P_0)=0$, and $\Omega$ being the differential \eqref{Omega} evaluated at ramification points according to \eqref{evaluation}, we have the following identity:
\begin{multline}
\label{s}
\frac{W(\pum, P_0)\Omega(P_0)}{\Omega(\pum)} + \sum_{i=1}^g\frac{W(\pum, \pxi)\Omega(\pxi)}{\Omega(\pum)}  + \sum_{\substack{j=1\\j\neq m}}^g \frac{W(\pum, \puj)\Omega(\puj)}{\Omega(\pum)}
\\
=\left( \sum_{i=1}^g\frac{\partial u_m}{\partial x_i}-1\right)\left(\frac{\prod_{\substack{i=1\\i\neq m}}^g(u_m-u_i)}{\prod_{i=1}^g(-u_i)} +\sum_{\substack{j=1\\j\neq m}}^g\frac{u_m \prod_{i\neq m,j}(u_m-u_i)}{u_j\prod_{i\neq j}(u_j-u_i)}\right)
\\
- \sum_{j=1}^g \frac{\prod_{\substack{i=1\\i\neq m}}(u_m-u_i)}{\prod_{i=1}^g(x_j-u_i)} \frac{\partial u_m}{\partial x_j}
- \sum_{\substack{j=1\\j\neq m}}^g\sum_{i=1}^g\frac{ (x_i-u_m) \prod_{s\neq m,j}(u_m-u_s) }{(x_i-u_j)\prod_{s\neq j}(u_j-u_s)}\frac{\partial u_m}{\partial x_i} - I_m^{u_m}.
\end{multline}
\end{lemma}
{\it Proof.} Let us denote by $S$ the quantity from the lemma which we seek to express:
\begin{equation*}
S:=\frac{W(\pum, P_0)\Omega(P_0)}{\Omega(\pum)} + \sum_{i=1}^g\frac{W(\pum, \pxi)\Omega(\pxi)}{\Omega(\pum)}  + \sum_{\substack{j=1\\j\neq m}}^g \frac{W(\pum, \puj)\Omega(\puj)}{\Omega(\pum)}
\end{equation*}
and note that $S\Omega(\pum)$ can be obtained through the sum of residues of the differential $\frac{\Omega(P)W(P, \pum)}{d\lambda(P)}$ on a compact Riemann surface from the family \eqref{hyper}:
\begin{equation*}
0=\sum_{P\in \mathcal C_\x}\underset{P}{\rm res}\frac{\Omega(P)W(P, \pum)}{d\lambda(P)}=\frac{1}{2}S\Omega(\pum) +
\underset{P=\pum}{\rm res}\frac{\Omega(P)W(P, \pum)}{d\lambda(P)}\,.
\end{equation*}
Now rewriting $W(P, \pum)$ as in \eqref{Wppum}, we have for the last residue:
\begin{equation*}
\underset{P=\pum}{\rm res}\frac{\Omega(P)W(P, \pum)}{d\lambda(P)}= \underset{P=\pum}{\rm res}\frac{\Omega(P)\phi(P)}{\phi(\pum)(\lambda(P)-u_m)d\lambda(P)} + \frac{1}{2}\Omega(\pum) I_m^{u_m}\,.
\end{equation*}
The sum of residues of the differential $\frac{\Omega(P)\phi(P)}{(\lambda(P)-u_m)d\lambda(P)}$ vanishes and thus we have
\begin{multline*}
0=\sum_{P\in\mathcal C_\x} \underset{P}{\rm res}\frac{\Omega(P)\phi(P)}{(\lambda(P)-u_m)d\lambda(P)}=
\underset{P=\pum}{\rm res}\frac{\Omega(P)\phi(P)}{(\lambda(P)-u_m)d\lambda(P)}
\\
+ \frac{1}{2}\left( \frac{\Omega(P_0)\phi(P_0)}{-u_m} + \sum_{j=1}^g \frac{\Omega(\pxj)\phi(\pxj)}{x_j-u_m} + \sum_{\substack{j=1\\j\neq m}}^g\frac{\Omega(\puj)\phi(\puj)}{u_j-u_m}\right)\,.
\end{multline*}
Putting these results together, we obtain
\begin{equation*}
S= \frac{1}{\Omega(\pum)\phi(\pum)}\left( \frac{\Omega(P_0)\phi(P_0)}{-u_m} + \sum_{j=1}^g \frac{\Omega(\pxj)\phi(\pxj)}{x_j-u_m} + \sum_{\substack{j=1\\j\neq m}}^g\frac{\Omega(\puj)\phi(\puj)}{u_j-u_m}\right) - I_m^{u_m}\,.
\end{equation*}
Let us rewrite the terms with $\Omega(\pxi)$ and $\Omega(P_0)$ with the help of differentials $v_j$ making use of their expression from \eqref{v}. In the last term, we express $\Omega(\puj)$ as in \eqref{Oum} using $j$ instead of $m$ in \eqref{Oum}. This, together with expression \eqref{derivative} for derivatives $\frac{\partial u_m}{\partial x_j} ,$ yields
\begin{multline*}
S=\frac{\Omega(P_0)v_m(P_0)}{\Omega(\pum) }\frac{\prod_{\substack{i=1\\i\neq m}}^g(u_m-u_i)}{\prod_{i=1}^g(-u_i)} -\prod_{\substack{i=1\\i\neq m}}(u_m-u_i) \sum_{j=1}^g \frac{1}{\prod_{i=1}^g(x_j-u_i)} \frac{\partial u_m}{\partial x_j}
\\
- \frac{1}{\phi(\pum)\Omega(\pum)}\sum_{\substack{j=1\\j\neq m}}^g\frac{\phi(\puj)(\Omega(P_0)v_j(P_0) + \sum_{i=1}^g\Omega(\pxi)v_j(\pxi))}{u_j-u_m} - I_m^{u_m}.
\end{multline*}
Rewriting the expressions in terms of differentials $v_j$ \eqref{v} once again and expressing $\frac{\Omega(P_0)v_m(P_0)}{\Omega(\pum) }$ via derivatives of $u_m$ due to \eqref{useful}, we prove the lemma.
$\Box$

\begin{lemma}
\label{lemma_t}
Let $m$ and $k$ be fixed integers, $1\leqslant m\leqslant g$ and let $I_j^{x_k}$ for $j=1, \dots, g$ be the normalization constants from \eqref{Wppxk}.
Assuming notation of Lemma \ref{lemma_s}, we have the following identity:
\begin{multline}
\label{t}
\frac{\partial u_m}{\partial x_k}\left( \frac{W(\pxk, P_0)\Omega(P_0)}{\Omega(\pxk)}+\sum_{\substack{j=1\\j\neq k}}^g \frac{W(\pxk, \pxj)\Omega(\pxj)}{\Omega(\pxk)}   +\sum_{j=1}^g \frac{W(\pxk, \puj)\Omega(\puj)}{\Omega(\pxk)}  \right)=
\\
\left( \sum_{i=1}^g\frac{\partial u_m}{\partial x_i}-1\right)\left(\frac{\prod_{i\neq m} (x_k-u_i)}{x_k\prod_{i\neq m} (-u_i)}
- \sum_{j=1}^g\frac{ u_m}{u_j} \frac{\prod_{i\neq m, j} (x_k-u_i)}{ \prod_{i\neq j} (u_j-u_i)} \right)
- \frac{\partial u_m}{\partial x_k}\sum_{j=1}^g I_j^{x_k}v_j(\pxk)
\\
+\sum_{\substack{j=1\\j\neq k}}^g\frac{ \prod_{i\neq m} (x_k-u_i)}{ \prod_{i\neq m} (x_j-u_i)}\frac{1}{(x_j-x_k)}\frac{\partial u_m}{\partial x_j}
%\\
- \sum_{i,j=1}^g \frac{\partial u_m}{\partial x_i}\frac{ x_i-u_m}{(u_j-x_k) (x_i-u_j)} \frac{ \prod_{s\neq m} (x_k-u_s)}{ \prod_{s\neq j} (u_j-u_s)}.
\end{multline}
\end{lemma}
{\it Proof.} Let us denote by $T$ the quantity in the parentheses in the left hand side:
\begin{equation*}
T:= \frac{W(\pxk, P_0)\Omega(P_0)}{\Omega(\pxk)}+\sum_{\substack{j=1\\j\neq k}}^g \frac{W(\pxk, \pxj)\Omega(\pxj)}{\Omega(\pxk)}   +\sum_{j=1}^g \frac{W(\pxk, \puj)\Omega(\puj)}{\Omega(\pxk)}.
\end{equation*}
Then note that $T\Omega(\pxk)$ can be expressed through the sum of residues of the differential $\frac{\Omega(P)W(P, \pxk)}{d\lambda(P)}$ on a compact Riemann surface corresponding to a curve from the family \eqref{hyper}:
\begin{equation*}
0=\sum_{P\in \mathcal C_\x}\underset{P}{\rm res}\frac{\Omega(P)W(P, \pxk)}{d\lambda(P)}=\frac{1}{2}T\Omega(\pxk) +
\underset{P=\pxk}{\rm res}\frac{\Omega(P)W(P, \pxk)}{d\lambda(P)}\,.
\end{equation*}
Representing $W(P, \pxk)$ in coordinates as in \eqref{Wppxk}, we have
\begin{equation*}
\underset{P=\pxk}{\rm res}\frac{\Omega(P)W(P, \pxk)}{d\lambda(P)}=\frac{1}{\phi(\pxk)}\underset{P=\pxk}{\rm res}\frac{\Omega(P)\phi(P)}{(\lambda(P)-x_k)d\lambda(P)} + \frac{1}{2} \Omega(\pxk) \sum_{j=1}^g I_j^{x_k}v_j(\pxk)\,.
\end{equation*}
The last residue can be obtained through the sum of residues of the differential $\frac{\Omega(P)\phi(P)}{(\lambda(P)-x_k)d\lambda(P)}$ vanishing on a compact Riemann surface. This differential has simple poles at all finite ramification points except $\pxk$ where it has a pole of order three. Thus we have
\begin{multline*}
0=\sum_{P\in \mathcal C_\x}\underset{P}{\rm res}\frac{\Omega(P)\phi(P)}{(\lambda(P)-x_k)d\lambda(P)}=\underset{P=\pxk}{\rm res}\frac{\Omega(P)\phi(P)}{(\lambda(P)-x_k)d\lambda(P)}+\frac{1}{2}\frac{\Omega(P_0)\phi(P_0)}{-x_k}
\\
+\frac{1}{2}\sum_{\substack{j=1\\j\neq k}}^g\frac{\Omega(\pxj)\phi(\pxj)}{x_j-x_k}+ \frac{1}{2}\sum_{j=1}^g\frac{\Omega(\puj)\phi(\puj)}{u_j-x_k}.
\end{multline*}
Putting these results together, we obtain
\begin{multline*}
T
%=- \frac{2}{\Omega(\pxk)}\underset{P=\pxk}{\rm res}\frac{\Omega(P)W(P, \pxk)}{d\lambda(P)}
%\\
=- \frac{\Omega(P_0)\phi(P_0)}{x_k\Omega(\pxk)\phi(\pxk)}
+\sum_{\substack{j=1\\j\neq k}}^g\frac{\Omega(\pxj)\phi(\pxj)}{\Omega(\pxk)\phi(\pxk)(x_j-x_k)}+ \sum_{j=1}^g\frac{\Omega(\puj)\phi(\puj)}{\Omega(\pxk)\phi(\pxk)(u_j-x_k)} - \sum_{j=1}^g I_j^{x_k}v_j(\pxk).
\end{multline*}
Let us now multiply $T$ by $\frac{\partial u_m}{\partial x_k}=-\frac{v_m(\pxk)\Omega(\pxk)}{\Omega(\pum)}$ to express the left hand side of the equation of the lemma. This leads to derivatives $\frac{\partial u_m}{\partial x_j}$ appear in one of the sums in the above expression for $T$. At the same time, we replace $\Omega(\puj)$ from \eqref{Oum} with $m=j.$ This yields
\begin{multline*}
T\frac{\partial u_m}{\partial x_k}
= \frac{v_m(\pxk)}{x_k\phi(\pxk)}\frac{\Omega(P_0)\phi(P_0)}{\Omega(\pum)}
+\sum_{\substack{j=1\\j\neq k}}^g\frac{v_m(\pxk)}{\phi(\pxk)(x_j-x_k)}\frac{\phi(\pxj)}{v_m(\pxj)}\frac{\partial u_m}{\partial x_j}
\\
+ \sum_{j=1}^g\frac{v_m(\pxk)\phi(\puj)}{\phi(\pxk)}\frac{(\Omega(P_0)v_j(P_0) + \sum_{i=1}^g\Omega(\pxi)v_j(\pxi))}{\Omega(\pum)(u_j-x_k)} - \frac{\partial u_m}{\partial x_k}\sum_{j=1}^g I_j^{x_k}v_j(\pxk).
\end{multline*}
Using explicit expressions \eqref{v} for differentials $v_j$, grouping the terms with $\frac{\Omega(P_0)v_m(P_0)}{\Omega(\pum)}$ and replacing this quantity by $\sum_{i=1}^g\frac{\partial u_m}{\partial x_i}-1$ due to \eqref{useful}, we prove the lemma.
$\Box$

\subsection{System of equations for the branch points under isoperiodic deformations}
\label{sect_main}

Now we are in  position to derive the system of equations for the functions $u_1, \dots, u_g$ depending on the variables $x_1, \dots, x_g$ which needs to be satisfied to provide a deformation for a Hill curve.
\begin{theorem}
\label{thm_main}
 Let $\Sx$ be a connected domain as discussed in Section \ref{sect_iso}, that is such that for $\x=(x_1,\dots, x_g)\in\Sx$ one can choose a canonical homology basis $\{a_j, b_j\}_{j=1}^g$ on all Riemann surfaces corresponding to the curves from the family \eqref{hyperx} in a consistent way meaning that $\lambda(a_j)$ and $\lambda(b_j)$ are constant for all $j$ and all $\x\in\Sx.$
Let $u_1(\x), \dots, u_g(\x)$ with $\x\in\Sx$ be such that the compact Riemann surfaces associated with curves of subfamily
\begin{equation}
\label{subfamily2}
\mathcal C_\x = \{(\l,\mu)\in\mathbb C^2 \;|\; \mu^2=\l\prod_{j=1}^g(\l-u_j(\x))\prod_{j=1}^g(\l-x_j)\}
\end{equation}
 of the family \eqref{hyper}
are such that the differential of the second kind $\Omega$ defined on each of these surfaces by \eqref{Omega} has $a$- and $b$-periods that are constant for all $\x\in\Sx$.
%along the subfamily \eqref{subfamily1}.
Then the functions $u_1(\x), \dots, u_g(\x)$ satisfy the following system of equations. For $1\leqslant k, m, n\leqslant g$ and $n\neq k$
\begingroup
\allowdisplaybreaks
\begin{multline}
\label{umxnxk}
\frac{\partial^2 u_m}{\partial x_k\partial x_n}
=\frac{1}{2}\frac{\partial u_m}{\partial x_k} \left( \frac{1}{x_k-x_n}+ \frac{1}{x_n-u_m} \right)
+\frac{1}{2}\frac{\partial u_m}{\partial x_n}\left( \frac{1}{x_n-x_k} + \frac{1}{x_k-u_m}\right)
\\
+ \frac{1}{2}\frac{\partial u_m}{\partial x_k}\frac{\partial u_m}{\partial x_n} \left(\frac{1}{u_m}+\sum_{\substack{i=1\\i\neq k,n}}^g \frac{1}{u_m-x_i}-\sum_{\substack{i=1\\ i\neq m}}^g \frac{2}{u_m-u_i} \right)
\\
+\frac{1}{4}\frac{\partial u_m}{\partial x_k}\sum_{\substack{j=1\\j\neq m}}^g \left(\frac{1}{u_m-u_j}-\frac{1}{x_k-u_j} \right) \frac{\partial u_j}{\partial x_n}
+\frac{1}{4}\frac{\partial u_m}{\partial x_n}\sum_{\substack{j=1\\j\neq m}}^g \left(\frac{1}{u_m-u_j}-\frac{1}{x_n-u_j} \right) \frac{\partial u_j}{\partial x_k}
\\
- \frac{1}{2}\frac{\partial u_m}{\partial x_k}\frac{\partial u_m}{\partial x_n}\left( \sum_{i=1}^g\frac{\partial u_m}{\partial x_i}-1\right)\prod_{\substack{i=1\\i\neq m}}^g\frac{u_i-u_m}{u_i}\left(\frac{1}{u_m} -\sum_{\substack{j=1\\j\neq m}}^g\frac{1}{u_m-u_j}\prod_{\substack{s=1\\s\neq j}}^g\frac{u_s}{u_s-u_j}\right)
\\
-\frac{1}{2}\frac{\partial u_m}{\partial x_k} \frac{\partial u_m}{\partial x_n} \left(\sum_{j=1}^g\frac{1}{x_j-u_m} \prod_{\substack{i=1\\i\neq m}}^g\frac{u_m-u_i}{x_j-u_i} \frac{\partial u_m}{\partial x_j}
+\sum_{j=1}^g \sum_{\substack{i=1\\i\neq m}}^g\frac{ x_j-u_m  }{(x_j-u_i)(u_m-u_i)}\frac{\prod_{s\neq m}(u_m-u_s)}{\prod_{s\neq i}(u_i-u_s)}\frac{\partial u_m}{\partial x_j}  \right)
\end{multline}
\endgroup
and
\begingroup
\allowdisplaybreaks
\begin{multline}
\label{umxkxk}
\frac{\partial^2 u_m}{\partial x_k^2}
=\frac{1}{2}\frac{\partial u_m}{\partial x_k} \left(
- \frac{1}{x_k}-\sum_{\substack{j=1\\j\neq k}}^g\frac{1}{x_k-x_j} + \sum_{\substack{j=1\\j\neq m}}^g \frac{2}{x_k-u_j} +\frac{1}{x_k-u_m}
\right)
 \\
  +\frac{1}{2}\left(\frac{\partial u_m}{\partial x_k} \right)^2\left(  \frac{1}{u_m}  + \sum_{\substack{j=1\\j\neq k}}^g \frac{1}{u_m-x_j}
 - \sum_{\substack{j=1\\j\neq m}}^g \frac{2}{u_m-u_j}+ \frac{1}{x_k-u_m}
 \right)
 \\
  +\frac{1}{2}\frac{\partial u_m}{\partial x_k}  \sum_{\substack{j=1\\j\neq m}}^g \left(\frac{1}{u_m-u_j}-\frac{1}{x_k-u_j} \right) \frac{\partial u_j}{\partial x_k}
 %%%%%%
 %+\frac{1}{2}\frac{\partial u_m}{\partial x_k}  \sum_{\substack{j=1\\j\neq m}}^g \frac{x_k-u_m}{(u_m-u_j)(x_k-u_j)}  \frac{\partial u_j}{\partial x_k}
 \\
-\frac{1}{2}\left( \sum_{i=1}^g\frac{\partial u_m}{\partial x_i}-1\right)\prod_{\substack{i=1\\i\neq m}}^g\frac{u_i-x_k}{u_i} \left(\frac{1}{x_k}
- \sum_{j=1}^g\frac{1}{x_k-u_j} \prod_{\substack{s=1\\s\neq j}}^g\frac{u_s}{u_s-u_j} \right)
\\
-\frac{1}{2} \left( \frac{\partial u_m}{\partial x_k}\right)^2\left( \sum_{i=1}^g\frac{\partial u_m}{\partial x_i}-1\right)\prod_{\substack{i=1\\i\neq m}}^g\frac{u_i-u_m}{u_i}\left(\frac{1}{u_m} -\sum_{\substack{j=1\\j\neq m}}^g\frac{1}{u_m-u_j}\prod_{\substack{s=1\\s\neq j}}^g\frac{u_s}{u_s-u_j}\right)
\\
-\frac{1}{2}\sum_{\substack{j=1\\j\neq k}}^g\frac{1}{x_j-x_k}\prod_{\substack{i=1\\i\neq m}}^g\frac{  x_k-u_i}{x_j-u_i}\frac{\partial u_m}{\partial x_j}
%\
-\frac{1}{2}\sum_{i,j=1}^g \frac{ x_j-u_m}{ (x_j-u_i)(x_k-u_m)} \prod_{\substack{s=1\\s\neq i}}^g\frac{x_k-u_s}{u_i-u_s} \frac{\partial u_m}{\partial x_j}
\\
-\frac{1}{2}\left( \frac{\partial u_m}{\partial x_k}\right)^2 \left(\sum_{j=1}^g \frac{1}{x_j-u_m}\prod_{\substack{i=1\\i\neq m}}^g\frac{u_m-u_i}{x_j-u_i} \frac{\partial u_m}{\partial x_j}
+ \sum_{j=1}^g\sum_{\substack{i=1\\i\neq m}}^g\frac{ x_j-u_m  }{(x_j-u_i)(u_m-u_i)}\frac{\prod_{s\neq m}(u_m-u_s)}{\prod_{s\neq i}(u_i-u_s)}\frac{\partial u_m}{\partial x_j}  \right).
\end{multline}
\endgroup
\end{theorem}
{\it Proof.} We subdivide the proof into two parts, proving the two equations claimed in the theorem in the two parts, respectively.

\begin{paragraph}{Part 1: proof of equation \eqref{umxnxk}.}Let us differentiate expression \eqref{derivative} for $\frac{\partial u_m}{\partial x_k}$ with respect to $x_n$ with $n\neq k:$
\begin{equation}
\label{der1-initial}
\frac{\partial^2 u_m}{\partial x_k\partial x_n}=\frac{\partial}{\partial x_n} \left\{  -  \frac{ v_{m}(\pxk)\Omega(\pxk)}{ \Omega(\pum)}   \right\}
=\frac{\partial u_m}{\partial x_k} \left( \frac{\frac{\partial v_m(\pxk)}{\partial x_n}}{v_m(\pxk)} + \frac{\frac{\partial \Omega(\pxk)}{\partial x_n}}{\Omega(\pxk)} -  \frac{\frac{\partial \Omega(\pum)}{\partial x_n}}{\Omega(\pum)}    \right).
\end{equation}
Let us write these three terms one by one. Note that derivatives with respect to $x_j$ contain derivatives with respect to the dependent branch points $u_m$.  Differentiation of $v_m(\pxk)$ can be done using explicit expressions \eqref{v} and \eqref{phi}, \eqref{phi-eval}:
\begin{multline}
\label{term1}
\frac{1}{v_m(\pxk)}\frac{\partial v_m(\pxk)}{\partial x_n} =\frac{1}{2}\left( \frac{1}{x_k-x_n}- \frac{1}{u_m-x_n}  -\sum_{\substack{j=1\\j\neq m}}^g\frac{1}{x_k-u_j}\frac{\partial u_j}{\partial x_n}
\right.
\\
\left.
+ \left(\frac{1}{u_m}+\sum_{\substack{i=1\\i\neq k}}^g \frac{1}{u_m-x_i}\right)\frac{\partial u_m}{\partial x_n}
- \sum_{\substack{j=1\\j\neq m}}^g \frac{1}{u_m-u_j} \left( \frac{\partial u_m}{\partial x_n}-\frac{\partial u_j}{\partial x_n}\right) \right).
\end{multline}
The Rauch formulas \eqref{Rauch-Omegalj} and the chain rule, imply the following expression for the derivative of $\Omega(\pxk):$
\begin{multline*}
\frac{1}{\Omega(\pxk)}\frac{\partial \Omega(\pxk)}{\partial x_n} = \frac{W(\pxk, \pxn)\Omega(\pxn)}{2\Omega(\pxk)} +\sum_{j=1}^g \frac{W(\pxk, \puj)\Omega(\puj)}{2\Omega(\pxk)}  \frac{\partial u_j}{\partial x_n}
\\
=\frac{\Omega(\pxn)}{2\Omega(\pxk)}\left(W(\pxk, \pxn) -\sum_{j=1}^g W(\pxk, \puj)v_j(\pxn)\right),
\end{multline*}
where in the second equality we rewrote the derivative $\frac{\partial u_j}{\partial x_n}$ according to \eqref{derivative}. Applying now Lemma \ref{lemma_T2}, we obtain
\begin{equation*}
\frac{1}{\Omega(\pxk)}\frac{\partial \Omega(\pxk)}{\partial x_n}
=\frac{\Omega(\pxn)}{2\Omega(\pxk)}\frac{\phi(\pxn)}{\phi(\pxk)} \frac{1}{x_n-x_k} \frac{\prod_{i=1}^g(x_n-u_i)}{\prod_{i=1}^g(x_k-u_i)},
\end{equation*}
and thus, using again \eqref{derivative} and \eqref{v}, we have
\begin{equation}
\label{term2}
\frac{\partial u_m}{\partial x_k}\frac{\frac{\partial \Omega(\pxk)}{\partial x_n} }{\Omega(\pxk)}
=-\frac{\Omega(\pxn)v_m(\pxn)}{2 \Omega(\pum)}  \frac{x_n-u_m}{(x_n-x_k)(x_k-u_m)}=\frac{1}{2}\frac{\partial u_m}{\partial x_n}\left( \frac{1}{x_n-x_k} + \frac{1}{x_k-u_m}\right).
\end{equation}

To differentiate $\Omega(\pum)$ we Corollary \ref{cor_epsilon} with $\{\l_1, \dots, \l_{2g+1}\}=\{0,x_1, \dots, x_g, u_1, \dots, u_g\}$ to obtain the derivative of $\Omega(\pum)$ with respect to $u_m:$
%need to use Lemma \eqref{lemma_epsilon} to express the derivative with respect to $u_m:$
%
\begin{multline*}
\frac{1}{\Omega(\pum)}\frac{\partial \Omega(\pum)}{\partial x_n}
=\frac{1}{\Omega(\pum)}\left( \frac{\partial^{\rm Rauch} \Omega(\pum)}{\partial x_n} + \sum_{\substack{j=1\\ j\neq m}}^g \frac{\partial^{\rm Rauch} \Omega(\pum)}{\partial u_j} \frac{\partial u_j}{\partial x_n} +\frac{\partial^{\rm Rauch} \Omega(\pum)}{\partial u_m}\frac{\partial u_m}{\partial x_n}\right)
\\
= \frac{W(\pum, \pxn)\Omega(\pxn)}{2\Omega(\pum)} + \frac{1}{2}\sum_{\substack{j=1\\j\neq m}}^g \frac{W(\pum, \puj)\Omega(\puj)}{\Omega(\pum)}\frac{\partial u_j}{\partial x_n}
\\
-\frac{1}{2}\left( \frac{W(\pum, P_0)\Omega(P_0)}{\Omega(\pum)} + \sum_{i=1}^g\frac{W(\pum, \pxi)\Omega(\pxi)}{\Omega(\pum)}  + \sum_{\substack{j=1\\j\neq m}}^g \frac{W(\pum, \puj)\Omega(\puj)}{\Omega(\pum)} \right) \frac{\partial u_m}{\partial x_n}.
\end{multline*}
After rewriting derivative $\frac{\partial u_j}{\partial x_n}$  in terms of $\Omega$ and $v_j$ as in \eqref{derivative}, we can apply Lemma \ref{lemma_T2} in the first line. In the second line, Lemma \ref{lemma_s} can be applied. This yields
\begin{multline}
\label{term3}
\frac{1}{\Omega(\pum)}\frac{\partial \Omega(\pum)}{\partial x_n} = \frac{\Omega(\pxn)v_m(\pxn)}{2\Omega(\pum)}\left(\frac{1}{x_n-u_m}- \sum_{\substack{i=1\\ i\neq m}}^g \frac{1}{u_m-u_i} + I_m^{u_m} \right)
\\
-\frac{1}{2} \left( \sum_{i=1}^g\frac{\partial u_m}{\partial x_i}-1\right)\left(\frac{\prod_{\substack{i=1\\i\neq m}}^g(u_m-u_i)}{\prod_{i=1}^g(-u_i)} +\sum_{\substack{j=1\\j\neq m}}^g\frac{u_m \prod_{i\neq m,j}(u_m-u_i)}{u_j\prod_{i\neq j}(u_j-u_i)}\right) \frac{\partial u_m}{\partial x_n}
\\
+\frac{1}{2} \left(\sum_{j=1}^g \frac{\prod_{\substack{i=1\\i\neq m}}(u_m-u_i)}{\prod_{i=1}^g(x_j-u_i)} \frac{\partial u_m}{\partial x_j}
+ \sum_{\substack{j=1\\j\neq m}}^g\sum_{i=1}^g\frac{ (x_i-u_m) \prod_{s\neq m,j}(u_m-u_s) }{(x_i-u_j)\prod_{s\neq j}(u_j-u_s)}\frac{\partial u_m}{\partial x_i} + I_m^{u_m} \right) \frac{\partial u_m}{\partial x_n}.
\end{multline}
Note that the terms containing $I_m^{u_m}$ cancel out due to \eqref{derivative}, and the remaining terms are rational in $u_j$ and their derivatives, also due to \eqref{derivative} as for the very first factor in the right hand side we have $\frac{\Omega(\pxn)v_m(\pxn)}{\Omega(\pum)}=-\frac{\partial u_m}{\partial x_n}$.

Let us now use \eqref{term1}, \eqref{term2}, \eqref{term3} in expression \eqref{der1-initial} for the derivative $\frac{\partial^2 u_m}{\partial x_k\partial x_n}.$ This yields a right hand side rational in $u_j$ and their derivatives.
To obtain equation \eqref{umxnxk} claimed in the theorem, it remains to symmetrize the following term. We show that this term is symmetric with respect to exchanging $k$ and $n$ by rewriting the derivatives in terms of $\Omega$ and $v_m$ as in \eqref{derivative}:
\begin{multline*}
\frac{\partial u_m}{\partial x_k}\sum_{\substack{j=1\\j\neq m}}^g \left(\frac{1}{u_m-u_j}-\frac{1}{x_k-u_j} \right) \frac{\partial u_j}{\partial x_n}
=
\frac{\Omega(\pxk)\Omega(\pxn)}{\Omega(\pum)}\sum_{\substack{j=1\\j\neq m}}^g \left(\frac{1}{u_m-u_j}-\frac{1}{x_k-u_j} \right) \frac{v_j(\pxn)v_m(\pxk)}{\Omega(\puj)}
\\
=\frac{\Omega(\pxk)\Omega(\pxn)}{\Omega(\pum)\phi(\pum)}\sum_{\substack{j=1\\j\neq m}}^g \frac{\phi(\pxk)\phi(\pxn)\prod_{s=1}^g(x_k-u_s)\prod_{i=1}^g(x_n-u_i)}{\Omega(\puj)\phi(\puj)(u_m-u_j)(x_n-u_j)(x_k-u_j)\prod_{i\neq j}(u_j-u_i)\prod_{s\neq m}(u_m-u_s)}.
\end{multline*}
The last expression being symmetric in $n$ and $k$, we may symmetrize the left hand side and write
\begin{multline*}
\frac{1}{2}\frac{\partial u_m}{\partial x_k}\sum_{\substack{j=1\\j\neq m}}^g \left(\frac{1}{u_m-u_j}-\frac{1}{x_k-u_j} \right) \frac{\partial u_j}{\partial x_n}
\\
= \frac{1}{4}\frac{\partial u_m}{\partial x_k}\sum_{\substack{j=1\\j\neq m}}^g \left(\frac{1}{u_m-u_j}-\frac{1}{x_k-u_j} \right) \frac{\partial u_j}{\partial x_n}
+\frac{1}{4}\frac{\partial u_m}{\partial x_n}\sum_{\substack{j=1\\j\neq m}}^g \left(\frac{1}{u_m-u_j}-\frac{1}{x_n-u_j} \right) \frac{\partial u_j}{\partial x_k}.
\end{multline*}
This explains the third line in \eqref{umxnxk} and finishes the proof of the first part of the theorem.
\end{paragraph}
\begin{paragraph}{Part 2: proof of equation \eqref{umxkxk}.}
Now we proceed to differentiation of expression \eqref{derivative} for $\frac{\partial u_m}{\partial x_k}$ second time with respect to the same variable:
\begin{equation}
\label{der2-initial}
\frac{\partial^2 u_m}{\partial x_k^2}=\frac{\partial}{\partial x_k} \left\{  -  \frac{ v_{m}(\pxk)\Omega(\pxk)}{ \Omega(\pum)}   \right\}
=\frac{\partial u_m}{\partial x_k} \left( \frac{\frac{\partial v_m(\pxk)}{\partial x_k}}{v_m(\pxk)} + \frac{\frac{\partial \Omega(\pxk)}{\partial x_k}}{\Omega(\pxk)} -  \frac{\frac{\partial \Omega(\pum)}{\partial x_k}}{\Omega(\pum)}    \right).
\end{equation}
As in Part 1, the derivative of $v_m(\pxk)$ is obtained by differentiating the evaluation at $\pxk$ of explicit expression \eqref{v} for $v_m$ using also expressions \eqref{phi-eval} for $\phi$ evaluated at ramification points:
\begin{multline}
\label{term1-2}
\frac{1}{v_m(\pxk)}\frac{\partial v_m(\pxk)}{\partial x_k}=\frac{1}{v_m(\pxk)}\frac{\partial}{\partial x_k}\left\{  \frac{\phi(\pxk)\prod_{i\neq m}(x_k-u_i)}{\phi(
\pum)\prod_{i\neq m}(u_m-u_i)} \right\}
= \frac{1}{2} \left(  \frac{1}{u_m}\frac{\partial u_m}{\partial x_k}  + \sum_{\substack{j=1\\j\neq k}}^g \frac{1}{u_m-x_j}\frac{\partial u_m}{\partial x_k}
\right.
\\
\left.
- \frac{1}{x_k}-\sum_{\substack{j=1\\j\neq k}}^g\frac{1}{x_k-x_j} + \sum_{\substack{j=1\\j\neq m}}^g \frac{1}{x_k-u_j} \left( 1-\frac{\partial u_j}{\partial x_k}\right) - \sum_{\substack{j=1\\j\neq m}}^g \frac{1}{u_m-u_j} \left(\frac{\partial u_m}{\partial x_k} -\frac{\partial u_j}{\partial x_k}\right)\right).
\end{multline}
To obtain the second term in \eqref{der2-initial} with the derivative of $\Omega(\pxk)$, we need to use the Rauch formulas \eqref{Rauch-Omegalj} and Corollary \ref{cor_epsilon} with $\{\l_1, \dots, \l_{2g+1}\}=\{0,x_1, \dots, x_g, u_1, \dots, u_g\}$ as in \eqref{hyper}:
\begin{multline}
\label{temp-part2}
\frac{1}{\Omega(\pxk)}\frac{\partial \Omega(\pxk)}{\partial x_k}
=\frac{1}{\Omega(\pxk)}\left( \frac{\partial^{\rm Rauch} \Omega(\pxk)}{\partial x_k} + \sum_{j=1}^g \frac{\partial^{\rm Rauch} \Omega(\pxk)}{\partial u_j} \frac{\partial u_j}{\partial x_k} \right)
= -\frac{W(\pxk, P_0)\Omega(P_0)}{2\Omega(\pxk)}
\\
   -\sum_{\substack{j=1\\j\neq k}}^g \frac{W(\pxk, \pxj)\Omega(\pxj)}{2\Omega(\pxk)}
-\sum_{j=1}^g \frac{W(\pxk, \puj)\Omega(\puj)}{2\Omega(\pxk)}
%\\
+\sum_{j=1}^g \frac{W(\pxk, \puj)\Omega(\puj)}{2\Omega(\pxk)} \frac{\partial u_j}{\partial x_k}.
\end{multline}
Let us replace in the last sum the derivative $\frac{\partial u_j}{\partial x_k}$ by its expression \eqref{derivative}. This allows us to write the last sum as $ -\frac{1}{2}\sum_{j=1}^g W(\pxk, \puj)v_j(\pxk)$ and to replace it by the right hand side of \eqref{T3} from Lemma \ref{lemma_T3}.
Let us at the same time multiply the whole expression by $\frac{\partial u_m}{\partial x_k}$ as in \eqref{der2-initial} in order to use Lemma \ref{lemma_t} for the remaining terms of \eqref{temp-part2}. After cancelling the terms containing the normalization constants $I_j^{x_k}$, this gives
\begin{multline}
\label{term2-2}
\frac{\partial u_m}{\partial x_k}\frac{\frac{\partial \Omega(\pxk)}{\partial x_k}}{\Omega(\pxk)} =-\frac{1}{2}\left( \sum_{i=1}^g\frac{\partial u_m}{\partial x_i}-1\right)\left(\frac{\prod_{i\neq m} (x_k-u_i)}{x_k\prod_{i\neq m} (-u_i)}
- \sum_{j=1}^g\frac{ u_m}{u_j} \frac{\prod_{i\neq m, j} (x_k-u_i)}{ \prod_{i\neq j} (u_j-u_i)} \right)
\\
-\frac{1}{2}\sum_{\substack{j=1\\j\neq k}}^g\frac{ \prod_{i\neq m} (x_k-u_i)}{ \prod_{i\neq m} (x_j-u_i)}\frac{1}{(x_j-x_k)}\frac{\partial u_m}{\partial x_j}
\\
-\frac{1}{2}\sum_{i,j=1}^g \frac{ x_i-u_m}{ x_i-u_j} \frac{ \prod_{s\neq m, j} (x_k-u_s)}{ \prod_{s\neq j} (u_j-u_s)} \frac{\partial u_m}{\partial x_i}
+\frac{1}{2}\frac{\partial u_m}{\partial x_k} \sum_{j=1}^g\frac{1}{x_k-u_j}.
\end{multline}
%
%To obtain the derivative in the third term of \eqref{der2-initial}, we also need to use Corollary \ref{cor_epsilon} with $\{\l_1, \dots, \l_{2g+1}\}=\{0,x_1, \dots, x_g, u_1, \dots, u_g\}$ to obtain the derivative of $\Omega(\pum)$ with respect to $u_m:$
The third term in \eqref{der2-initial}, the quantity $\frac{1}{\Omega(\pum)}\frac{\partial \Omega(\pum)}{\partial x_k},$ is already obtained in Part 1: we need to use the result in \eqref{term3} upon replacing the index $n$ by $k$.
%%
%\begin{multline*}
%\frac{1}{\Omega(\pum)}\frac{\partial \Omega(\pum)}{\partial x_k}
%=\frac{1}{\Omega(\pum)}\left( \frac{\partial^{\rm Rauch} \Omega(\pum)}{\partial x_k} + \sum_{\substack{j=1\\ j\neq m}}^g \frac{\partial^{\rm %Rauch} \Omega(\pum)}{\partial u_j} \frac{\partial u_j}{\partial x_k} +\frac{\partial^{\rm Rauch} \Omega(\pum)}{\partial u_m}\frac{\partial u_m}{\partial x_k}\right)
%\\
%= \frac{W(\pum, \pxk) \Omega(\pxk)}{2\Omega(\pum)} + \sum_{\substack{j=1\\ j\neq m}}^g \frac{W(\pum, \puj) \Omega(\puj)}{2\Omega(\pum)}\frac{\partial u_j}{\partial x_k}
%\\
%-\left( \frac{W(\pum, P_0)\Omega(P_0)}{2\Omega(\pum)}   +\sum_{j=1}^g \frac{W(\pum, \pxj)\Omega(\pxj)}{2\Omega(\pum)}
%+\sum_{\substack{j=1\\j\neq m}}^g \frac{W(\pum, \puj)\Omega(\puj)}{2\Omega(\pum)}\right)\frac{\partial u_m}{\partial x_k}
%\end{multline*}
%
It remains now to plug the results \eqref{term1-2}, \eqref{term2-2} and the appropriately modified \eqref{term3} into \eqref{der2-initial}. $\Box$
\end{paragraph}
\begin{example}
\label{examplegenus2}
If $g=2$, system \eqref{umxnxk}, \eqref{umxkxk} from Theorem \ref{thm_main} has the form:
\begin{multline}
\label{umx1x2}
%m=1, k=1, n=2
\frac{\partial^2 u_1}{\partial x_1\partial x_2}
=\frac{1}{2}\frac{\partial u_1}{\partial x_1} \left( \frac{1}{x_1-x_2}+ \frac{1}{x_2-u_1} \right)
+\frac{1}{2}\frac{\partial u_1}{\partial x_2}\left( \frac{1}{x_2-x_1} + \frac{1}{x_1-u_1}\right)
+ \frac{1}{2}\frac{\partial u_1}{\partial x_1}\frac{\partial u_1}{\partial x_2} \left(\frac{2}{u_1}- \frac{3}{u_2-u_1} \right)
\\
+\frac{1}{4}\frac{\partial u_1}{\partial x_1}\frac{\partial u_2}{\partial x_2}\left(\frac{1}{u_1-u_2}-\frac{1}{x_1-u_2} \right)
+\frac{1}{4}\frac{\partial u_1}{\partial x_2}\frac{\partial u_2}{\partial x_1} \left(\frac{1}{u_1-u_2}-\frac{1}{x_2-u_2} \right)
\\
- \frac{1}{2}\left(\frac{\partial u_1}{\partial x_1}\right)^2\frac{\partial u_1}{\partial x_2}\left( \frac{1}{u_1} +\frac{1}{x_1-u_1} \right)
- \frac{1}{2}\frac{\partial u_1}{\partial x_1}\left(\frac{\partial u_1}{\partial x_2}\right)^2\left( \frac{1}{u_1}  + \frac{1}{x_2-u_1} \right)
\end{multline}
and
\begin{multline}
\label{umx1x1}
\frac{\partial^2 u_1}{\partial x_1^2}
=\frac{1}{2}\left(\frac{1}{x_1}- \frac{1}{x_1-u_1}  \right)+\frac{1}{2}\frac{\partial u_1}{\partial x_1} \left(
- \frac{2}{x_1}-\frac{1}{x_1-x_2} +  \frac{1}{x_1-u_2} +\frac{1}{x_1-u_1}
\right)
-\frac{1}{2}\frac{\partial u_1}{\partial x_2}\left(\frac{1}{x_1}+\frac{1}{x_2-x_1}\right)
 \\
  +\frac{1}{2}\left(\frac{\partial u_1}{\partial x_1} \right)^2\left(  \frac{2}{u_1}  + \frac{1}{u_1-x_2}
 -  \frac{1}{u_1-u_2}+ \frac{1}{x_1-u_1} \right)
 %\\
  +\frac{1}{2}\frac{\partial u_1}{\partial x_1}  \frac{\partial u_2}{\partial x_1} \left(\frac{1}{u_1-u_2}-\frac{1}{x_1-u_2} \right)
 %%%%%%
\\
-\frac{1}{2} \left( \frac{\partial u_1}{\partial x_1}\right)^3\left(\frac{1}{u_1} + \frac{1}{x_1-u_1}\right)
-\frac{1}{2} \left( \frac{\partial u_1}{\partial x_1}\right)^2 \frac{\partial u_1}{\partial x_2}\left(\frac{1}{u_1} +\frac{1}{x_2-u_1}\right).
\end{multline}
Replacing $u_1$ by $u_2$ and vice versa in \eqref{umx1x2} and in \eqref{umx1x1} we obtain equations for $\frac{\partial^2 u_2}{\partial x_1\partial x_2}$ and for $\frac{\partial^2 u_2}{\partial x_1^2}$, respectively.
Replacing $x_1$ by $x_2$ and vice versa in \eqref{umx1x1} we obtain the equation for $\frac{\partial^2 u_1}{\partial x_2^2}$.
\end{example}
\begin{theorem}
\label{thm_converse-g}
Let $\x_0\in\mathbb C^g\setminus\{(x_1, \dots, x_g)\; |\; \exists\,i\neq j \text{  with }   x_{i}=x_j \}$ be  such that $x_{0j}\neq 0$ for $j=1, \dots, g$. Let $\{a_1, \dots, a_g,b_1, \dots, b_g\}$ be a canonical homology basis on $\mathcal C_{\x_0}$ such that the projections $u(a_j)$ and $u(b_j)$ with $j=1, \dots, g$ do not intersect a certain neighbourhood $\hat \X$ of $\x_0$.  Let $\Omega=\Omega(\x_0)$ be the meromorphic differential of the second kind defined by \eqref{Omega} on a compact Riemann surface $\mathcal C_{\x_0}$ of the hyperelliptic curve \eqref{hyperx} of genus $g$ with $\x=\x_0$.  Assume that $\x_0$ is such that $\Omega(\x_0; \puj)\neq 0$ for $j=1, \dots, g$. Then there exists a unique continuous family  $(\mathcal C_\x, \Omega(\x))$ providing an isoperiodic deformation of the pair $(\mathcal C_{\x_0}, \Omega(\x_0))$ for $\x$ varying in some neighbourhood $\X\subset\hat\X$ of $\x_0.$
\end{theorem}
{\it Proof.}
The defining condition for isoperiodic deformations of $(\mathcal C_{\x_0}, \Omega(\x_0))$ are given by  $\beta_k=\oint_{b_k}\Omega(x)={\rm const},$ $k=1, \dots, g$.  By the implicit function theorem, these relations define functions $u_1(\x), \dots, u_g(\x)$ in a neighbourhood of $\x=\x_0$ if the the Jacobian ${\rm det} J $ for a $g\times g$ matrix $J$ such that $J_{jk}= \frac{\partial \beta_k}{\partial u_j}$ is non-zero at $\x=\x_0$. From \eqref{betak} and \eqref{Omega} we obtain these derivatives using the Rauch formulas \eqref{RauchB}:
\begin{equation*}
J_{jk} = \pi\i\, W(P_\infty, \puj) \omega(\puj) + \pi\i\, \alpha\omega(\puj) \omega_k(\puj) = \pi\i\,\Omega(\puj)\omega_k(\puj).
\end{equation*}
Thus ${\rm det} J = \pi\i\,{\rm det}\hat J \prod_{m=1}^g\Omega(\pum)$ where by $\hat J$ we denote the matrix $\hat J_{jk} = \omega_k(\puj).$ The product of $\Omega(\pum)$ being non-zero at $\x=\x_0$, we only need to show that ${\rm det} \hat J\neq 0.$ Note that replacing the basis $\omega_1, \dots, \omega_g$ by the basis $\phi(P), \;\lambda(P)\phi(P), \dots, \lambda(P)^{g-1}\phi(P)$ \eqref{phi} in the space of holomorphic differentials does not affect the non-vanishing of ${\rm det}\hat J$. On the other hand, this change of basis transforms ${\rm det}\hat J$ into the Vandermonde determinant $V(u_1, \dots, u_g)$, which is non-zero for distinct $u_1, \dots, u_g$. This is fulfilled in a neighbourhood of $\x_0$. This proves the existence of local continuous isoperiodic deformation of $(\mathcal C_{\x_0}, \Omega(\x_0))$.

Now, Theorem \ref{thm_main} shows that for every continuous isoperiodic deformation, the functions $u_1(\x), \dots, u_g(\x)$ satisfy equations \eqref{umxnxk} and \eqref{umxkxk} for all $1\leqslant m,n,k\leqslant g$. In addition, due to Theorem \ref{thm_derivatives}, the condition $\oint_{b_k}\Omega(x)={\rm const}$ implies that the first derivatives of $\partial_{x_k}u_j(\x)$ are given by \eqref{derivative} for any $1\leqslant k,j\leqslant g$ and $\x=(x_1,\dots, x_g)$. Thus, for $\x_0=(x_1^0, \dots, x_g^0)$, the  initial values $u_j(\x_0)$ and derivatives $\partial_{x_k}u_j(\x)$ \eqref{derivative} for any $\x$ specify a unique solution of system \eqref{umxnxk}, \eqref{umxkxk} as follows. On the set of points $(x_1, x_2^0 \dots, x_g^0)$ for $x_1$ varying in some neighbourhood of $x_1^0$, the system reduces to a system of ordinary differential equations of second order with respect to the variable $x_1$, and thus the initial values $u_j(\x_0)$ and the derivatives $\partial_{x_1}u_j(\x_0)$ single out a unique solution $u_j(x_1, x_2^0 \dots, x_g^0)$ for $j=1, \dots, g$ and $x_1\in S_1$ for some open neighbourhood $S_1$ of $x_1^0$. Now, the obtained functions $u_j(x_1, x_2^0 \dots, x_g^0)$ give us initial values for the system of ordinary differential equations of second order with respect to the variable $x_2$ at any point of the set $\{(x_1, x_2^0 \dots, x_g^0)|\; x_1\in S_1\}$. Together with the values of derivatives $\partial_{x_2}u_j((x_1, x_2^0 \dots, x_g^0))$, they single out a unique solution $u_j(x_1, x_2, x_3^0 \dots, x_g^0)$ for $j=1, \dots, g$ and $(x_1,x_2)\in S_2$ for some open neighbourhood $S_2$ of $(x_1^0, x_2^0)$. Proceeding in this way, we obtain uniqueness of an isoperiodic deformation of $(\mathcal C_{\x_0}, \Omega(\x_0))$ locally in an open neighbourhood of $\x_0$.

$\Box$

\begin{corollary}
\label{cor_real}
Let the situation and notation be those of Theorem \ref{thm_converse-g}.
 Let the initial  hyperelliptic curve $\mathcal C_{\x_0}$ \eqref{hyperx} of genus $g$ with $\x=\x_0$  have all branch points $x_j$ and $u_j$ real and satisfy  the Hill condition \eqref{eq:hill}.
 Then, there exists  a canonical homology basis on $\mathcal C_{\x_0}$ such that there exists a unique continuous isoperiodic deformation of the pair $(\mathcal C_{\x_0}, \Omega(\x_0))$ for $\x$ varying in some real neighbourhood $\X\cap \mathbb R^g \subset\hat\X$ of $\x_0$ for which all the functions $u_j=u_j(x_1, \dots, x_g)$ are real-valued, and the obtained hyperelliptic curves satisfy the  Hill  condition \eqref{eq:hill} with $T$ being the same  for all curves.
\end{corollary}

{\it Proof.}
The reality condition follows from  \cite{Mei} and the uniqueness of  the isoperiodic deformation proven in Theorem \ref{thm_converse-g}. Theorem 3 of \cite{Mei} shows that Hill curves can alternatively be characterized through the spectrum of Schr\"odinger operators $-\frac{\partial^2}{\partial X^2} + v(X)$ with periodic potentials $v$, see also \cite{MM}. Using the notation of \cite{Mei} and  \cite{Nov}, the set of branch points of the curves \eqref{hyperx} is denoted by $\{E_0, E_1, \dots, E_{2g}\}\subset \mathbb R$ with $E_{2j-2}<E_{2j-1}<E_{2j}$ for $j=1, \dots, g$, where the intervals of the form $(E_{2j}, E_{2j+1})$ form the Bloch spectrum of the Schr\"odinger operator while the intervals of the form $[E_{2j-1},E_{2j}]$ are called {\it gaps}. Let us now choose the canonical homology basis on $\mathcal C_{\x_0}$ in such a way that the projections on the $\lambda$-sphere of the $a_j$ cycle goes once around the gap interval $[E_{2j-1},E_{2j}]$ for all $j=1, \dots g$. The projection of the cycle $b_j$ may be chosen to go once around the set of intervals $[E_0,E_1], \dots [E_{2j-2},E_{2j-1}]$.  For this choice, the cycles $a_j$ and $b_j$ satisfy relations $\rho a_j=a_j$, $\rho b_j=-b_j$  with $\rho$ being the antiholomorphic involution from Example \ref{ex:realhyp}.

 Theorem 3 of \cite{Mei} states that the set of branch points of the hyperelliptic curve \ref{hyperx} corresponds to a periodic Schr\"odinger operator if and only if equations (1.23), (1.24) of \cite{Mei} are satisfied by the differential $dp(\lambda)$ of the {\it  quasi-momentum} function $p(\lambda)$ as well as the asymptotic at the point at infinity given by $p(\lambda-E_0)=\sqrt{\lambda-E_0}$.  This asymptotic implies that $dp(\lambda) = -\frac{d\zeta_\infty}{\zeta_\infty^2}$ in terms of the local coordinate \eqref{coordinates}. With the above choice of $a$- and $b$-cycles, equations (1.24), (1.23) from \cite{Mei} imply that $\oint_{a_j}dp=0$ and $\oint_{b_j}dp=\frac{2\pi \i n_j}{T}$, respectively.  These conditions imply that $dp=-\Omega_0$ where $\Omega_0$ is the differential on $\mathcal C_{\x_0}$ given by \eqref{Omega} with $\alpha=0$, see also formula (29), p. 145 in \cite{Nov}. These conditions also imply that the corresponding surface is a Hill curve.

As explained in \cite{Mei}, Hill curves may also be characterized as solutions to the ``$g$-band extremal problem", which is a minimum deviation problem for holomorphic functions on the set of intervals $[E_{0}, E_{1}], \dots, [E_{2g-2}, E_{2g-1}],$  $[E_{2g}, \infty]$, where the endpoints of the intervals are also to be found. In Section 4 of \cite{Mei}, it is shown that one may fix $g+1$ real parameters consisting of $E_0$ (we set $E_0=0$ in this paper), $E_{2g}$ and one endpoint of each of the intervals $[E_{2j}, E_{2j+1}]$, for $j=1,\dots, g-1$ and find the other $g$ endpoints, that is $E_1$ and the second endpoint of each of the intervals $[E_{2j}, E_{2j+1}]$, for $j=1,\dots, g-1$, as real functions of the  $g+1$ parameters (one of which is here set equal to zero) in order for the set of intervals to carry a holomorphic function with minimal deviation.  For this discussion, see \cite{Mei} p. 843, paragraphs around formula (4.17).  In our case, we know that $T$ (denoted by $a$ in \cite{Mei}) is ``large enough'', because we start here from a Hill curve $\mathcal C_{\x_0}$ that determines the value of $T$. This result implies the existence of a family of Hill curves parametrized by $g$ real parameters (as we set here $E_0=0$) for which we may take $g$ branch points of the curves. By uniqueness proven in Theorem \ref{thm_converse-g}, these deformations coincide with isoperiodic deformations of the pair $(\mathcal C_{\x_0}, \Omega_0)$ restricted to the real line.
$\Box$

\section{Isoperiodic deformations of two-gap potentials}
\label{sec:twogap}

 Let us now show how Theorems \ref{thm_main} and \ref{thm_converse-g} as well as Corollary \ref{cor_real} imply the existence of isoperiodic deformations of finite-gap potentials.
According to the theory of finite-gap potentials, developed starting from mid 1970's by Novikov, Dubrovin, Lax, Its, Matveev, Mumford, McKean, Krichever, van Moerbeke, and others, see e.g. \cite{DMN}, all such potentials can be obtained as stationary solutions of the hierarchy of Korteweg - de Vries equations. For the case of two-gap potentials, following \cite{DMN}, \cite{Dub1981}, \cite{Nov}, let us consider a pair of linear differential operators
$\mathbb L, \mathbb A$ satisfying the Novikov equation
$$
[\mathbb L, \mathbb A]=0.
$$
We consider the operators
\begin{equation}
\label{L2}
\mathbb L=\frac{d^2}{dX^2}+v(X)
\end{equation}
and
$$
\mathbb A=16\frac{d^5}{dX^5}+20\left(v\frac{d^3}{dX^3} + \frac{d^3}{dX^3} v\right)+30v\frac{d}{dX}v-5\left(v''\frac{d}{dX}+\frac{d}{dX}v''\right)+c_1\left(4\frac{d^3}{dX^3}+3\Big(v\frac{d}{dX}+\frac{d}{dX}v\Big)\right)+c_2\frac{d}{dX}.
$$
The resulting Novikov equation for $v(x)$ has a well-known Lagrangian form with the Lagrangian
$$
\mathcal L(v, v', v'') = \frac{v''}{2}-\frac{5}{2}v''v^2+\frac{5}{2}v^4+c_1\left(\frac{v'^2}{2}+v^3\right)+c_2v^2+c_3v.
$$
It is equivalent to a Hamiltonian system (see \cite{Nov})  with the Hamiltonian
\begin{equation}
\label{H}
H(q_1, q_2, p_1, p_2)=p_1p_2 + V(q_1, q_2),
\end{equation}
where the potential is $V(q_1, q_2)=\frac{1}{8}(-5q_1^4-20q_1^2q_2+4c_2q_1^2-4q_2^2+c_3q_1)$. Here we assume that $c_1=0,$ which can be accomplished by a shift of $v$ by a constant and we set
$$
q_1=v, \quad q_2=v''-5v^2,  \quad p_1=q_1', \quad  p_2=q_1'.
$$
This Hamiltonian system has two first integrals, the energy $H$ and
$$
J=p_1^2+p_2^2(2q_2-c_2)+ 2q_1p_1p_2 +q_1^5+c_2q_1^3-4q_1q_2^2+2c_2q_1q_2+2c_2q_2.
$$
Thus, the system is completely integrable. Consider the polynomial
$$
\P_5(\lambda)=\lambda^5+\frac{c_2\lambda^3}{8}+\frac{c_3\lambda^2}{16}+\frac{1}{32}\left(H+2c_2^2\right)\lambda+\frac{J-c_2c_3}{2^8},
$$
and the associated genus two hyperelliptic curve
\begin{equation}\label{eq:gamma}
\Gamma=\{(\l,w)\in\mathbb C^2 \;|\; w^2=\P_5(\lambda)\}.
\end{equation}
Denote the zeros of $\P_5$ by $\lambda_j$, $j=1, \dots, 5$.
 The compactification of the ramified covering $\lambda: \Gamma\to\mathbb C$, $\lambda(\lambda, w)=\lambda$ is ramified at the point $P_{\infty}$ over $\lambda=\infty$.
Denote  the other ramification points by $P_j=(\lambda=\lambda_j, w=\P_5(z_j))$  and denote  the compact Riemann surface associated with the curve \eqref{eq:gamma} by $\Gamma$ as well. Let us fix a canonical homology basis $\{a_1, a_2; b_1,  b_2\}$  on $\Gamma$ as follows. Denote by $a_1$ the contour going around the points $P_{3}$ and $P_{2}$, by $a_2$ the contour around $P_{5}$ and $P_{4}$, by $b_1$ the contour going around $P_2$ and $P_1$, and by $b_2$ the contour going around the points $P_1, \; P_2,\; P_3$ and $P_4$. We may think of the contours $b_1, b_2$ as staying on one of the sheets of the ramified covering of the $\lambda$-sphere. The directions on the contours are chosen in a way to have the intersection indices
%$a_k\circ a_j=b_k\circ b_j=0$ and
$a_k\circ b_j=\delta_{kj}$.  Denote by $\omega_1$ and $\omega_2$ the Abelian differentials of the first kind  normalized with respect to the chosen homology basis.

 The integration \cite{Nov} of the Novikov equation and the Hamiltonian system with Hamiltonian $H$  \eqref{H}
 %take the form%get the Abel map form
%$$
%\dot{\gamma}_1=2\i\frac{\sqrt{\P_5(\gamma_1)}}{\gamma_1-\gamma_2}, \quad \dot\gamma_2=2\i\frac{\sqrt{\P_5(\gamma_2)}}{\gamma_2-\gamma_1},
%$$
%where $\gamma_1, \gamma_2$ are new coordinates at the common level set of the first integrals, defined through their symmetric functions
%$$
%\gamma_1 + \gamma_2=\frac{q_1}{2}, \qquad \gamma_1  \gamma_2=\frac{8q_1^2+q_2}{8}+ \frac{1}{2}\sum_{j<k}\lambda_j\lambda_k.
%$$
 results in the following solution
 \begin{equation}
 \label{v-potential}
 v(X)=2\frac{\partial^2}{\partial X^2}{\rm log}\,\theta(XU+z_0) + \text{const}.
 \end{equation}
 Here $U$ is, up to a constant factor, the vector of $\b$-periods of a differential of the second kind having a pole of order $2$ at $P_{\infty}$ with principal part $d{\sqrt{\lambda}}$, no other poles, and zero $a$-periods. In  our notation, this differential is $-\Omega_0$ where $\Omega_0$ is defined by \eqref{Omega} with $\alpha=0.$ More precisely,
 \begin{equation}
 \label{U}
 U=-\frac{1}{2\pi\i}\left(\oint_{b_1}\Omega_0, \oint_{b_2}\Omega_0\right)^T.
 \end{equation}
  Every solution $v(X)$ is codified by a two-dimensional vector
 \begin{equation}\label{eq:z0}
 z_0=-\mathcal A(Q_1+Q_2)-K,
 \end{equation}
 where $Q_1,\,Q_2$ is a non-special divisor on $\Gamma$ and  $K$ is the vector of Riemann constants and $ \mathcal A(P)$ is the Abel map, based at some point $P_0\in\Gamma$,
 $$
 \mathcal A(P)=\Big(\int_{P_0}^P\omega_1,\int_{P_0}^P\omega_2\Big)^T.
 $$
 %Here $Q_j=(\gamma_j,\sqrt{\P_5(\gamma_j)})$, $j=1, 2$.
 In the real case we assume $\lambda_5<\lambda_4<\dots<\lambda_1$.
 % and $\gamma_j$ are reals such that $\lambda_5< \gamma_1<\lambda_4$ and $\lambda_3< \gamma_2<\lambda_2$.
 The spectrum of $\mathbb L$ as an operator on $\mathcal L_2(-\infty, \infty)$ is the ray $(-\infty, \lambda_1]$ without two {\it gap intervals}  $(\lambda_{5}, \lambda_{4})$ and $(\lambda_{3}, \lambda_{2})$. This is why $v(X)$ is called a two-gap potential. The Riemann surface $\Gamma$ is called the spectral curve of the operator $\mathbb L$. It is natural to consider the eigenfunction of $\mathbb L$, denoted by $\psi$, such that
 \begin{equation*}
 \mathbb L\psi=\lambda\psi.
\end{equation*}
 This is a Baker-Akhiezer function, meromorphic on $\Gamma\setminus P_{\infty}$, with poles at $Q_1$ and $Q_2$, of the form
 \begin{equation}
  \label{eigenL}
 \psi(X, P)=\exp \Big(-X\int_{P_0}^P\Omega_0\Big)\frac{\theta(\mathcal A(P)+XU+z_0)\theta(z_0)}{\theta(\mathcal A(P)+z_0)\theta(XU+z_0)}.
 \end{equation}
 The group of periods with respect to $X$ of the derivative of ${\rm log}\,\psi$ coincides with the group of periods of the potential $v(X)$ \eqref{v}.
 %The Riemann surface $\Gamma$ is the spectral curve of the operator $\mathbb L$.
 An important part of modern theory of integrable systems is a study of isospectral deformations, that is deformations of $v$ which preserve $\Gamma$ and thus the spectrum of deformed operators $\mathbb L$. Here we  study {\it isoperiodic}  deformations, that is those that deform $v$ and $\Gamma$, while preserving $U$, the period of the deformed differential $\Omega_0$.  Here $U$ can be seen as the wavevector of solutions $\psi$. In particular, in periodic cases, preserving the wavevector may lead to preserving the period (or one of the periods) of $\psi$. Thus, our isoperiodic deformations are {\it wavevector-preserving} and in the periodic cases, they can be {\it period-preserving}.

 The analogous construction for hyperelliptic curves of genus $g\geqslant 1$,
 $$\Gamma_g= \{(\l,w)\in\mathbb C^2 \;|\;  w^2=\P_{2g+1}(\lambda)\},$$
 where $\lambda_{2g+1}<\dots<\lambda_1$ are the zeros of the polynomial $\P_{2g+1}$, generates $g$-gap potentials $v(X)$, see \cite{Nov}, \cite{Dub1981}.
  An obtained $g$-gap potential $v(X)$ is periodic with period $T$ if and only if the vector $U$ \eqref{U} satisfies
\begin{equation}
\label{eq:periodicity}
TU=(n_1e_1+\dots+ n_ge_g) + \mathbb B(m_1e_1+\dots+m_ge_g),
\end{equation}
with some integers $n_j, \;m_j$, $j=1,\dots, g$,   where $e_j$ is the $j$-th vector of the canonical basis of $\mathbb C^g$ and $\mathbb B$ is the Riemann matrix.  In particular, for $m_1=\ldots=m_g=0$, one recognizes the Hill curves condition \eqref{eq:hill}.
 \begin{proposition}
 \label{prop:2gap}
Let $v(X)$ be  a  real periodic $g$-gap potential corresponding to the spectral curve
 \begin{equation}
 \label{eq:gammareal}
 \Gamma= \{(\l,w)\in\mathbb C^2 \;|\;  w^2=\prod_{j=1}^{2g+1}(\lambda-\lambda_j)\}\qquad\text{with}\quad \lambda_{2g+1}<\lambda_{2g}<\ldots<\lambda_1\,.
\end{equation}
Assume that the $a$- and $b$-cycles of the spectral curve are selected as in Example \ref{ex:realhyp}.
 Denote  by $P_\infty$ the unique point at infinity of the compact Riemann surface $\Gamma$ of this hyperelliptic curve and let $\Omega_0$ be the differential on $\Gamma$ defined by \eqref{Omega} with $\alpha=0$.  Introduce $x_j=\lambda_{2j-1}-\lambda_{2g+1}$, $u_j=\lambda_{2j}-\lambda_{2g+1}$, $j=1,\dots, g$.   The pair $(\Gamma, \Omega_0)$ is subject to an isoperiodic deformation if and only if $u_1,\dots, u_g$ as functions of $x_1,\dots, x_g$ satisfy equations \eqref{umxnxk} and \eqref{umxkxk} of Theorem \ref{thm_main} and the initial conditions \eqref{derivative},  where $\Omega=\Omega_0$ and $v_j$ are defined by \eqref{v}. These deformations preserve reality  of $x_j$, $u_j$ and also preserve  the wavevector of $\psi(X)$ \eqref{eigenL}, of the eigenfunction of the operator $\mathbb L$ \eqref{L2}.  If the initial potential is periodic with the period $T>0$,    then the potentials corresponding to the deformed curves are also periodic with the  period $T$.  \end{proposition}
{\it Proof.}
The first part of the Proposition follows directly from Theorems  \ref{thm_derivatives},  \ref{thm_main}  and \ref{thm_converse-g}. The wavevector $U$ is related to the vector of $b$-periods of the differential $\Omega_0$ via equation \eqref{U} and it is further related to the normalized holomorphic differentials via the vector form of  \eqref{Omega0}:
\begin{equation}
\label{Omega0g}
U=\frac{1}{2\pi\i}\oint_b\Omega_0=\omega(P_\infty),
\end{equation}
where  $b=(b_1, \dots, b_g)$ and  $\omega(P_\infty)$ is the column consisting of the normalized differentials of the basis of holomorphic differentials, evaluated at $P_\infty$ according to \eqref{evaluation}. Due to the reality  of the branch points, $\Gamma$ is a real curve, see Example \ref{ex:realhyp}, and  thus $\omega_j(P_\infty)$ is real for all $j=1,\dots, g$.  Therefore, $U$ is a real vector, and it is preserved under the isoperiodic deformations of $(\Gamma, \Omega_0)$. These deformations preserve reality  of $x_j$, $u_j$ according to Corollary \ref{cor_real}.  Thus, if the initial potential \eqref{v-potential} is periodic with the period $T$, then $TU$ satisfies \eqref{eq:periodicity}. Since for a real curve with the given choice of canonical homology basis the Riemann matrix $\mathbb B$ is purely imaginary, then $TU$ satisfies \eqref{eq:periodicity}  with zero $m_j$'s. Therefore, $T$ remains to be a period under the isoperiodic deformations, since $U$ is real and remains unchanged under the deformations.
$\Box$
\begin{example}\label{ex:2gap}[Deformations of two-gap Lam\'e potentials] It was shown in \cite{In}  that Lam\'e potentials $v_g(x)=g(g+1)\wp(x)$ have $g$ gaps. In particular $v_2(x)=6\wp(x)$ is a two-gap potential.
In the real case, denoting by $e_1< e_2 < e_3$  the zeros of the polynomial $4\lambda^3-g_2\lambda-g_3$, following \cite{DMN}, the endpoints of the gaps are given by:
$$
\lambda_5=-3(e_2+e_3),\quad \lambda_4=-\sqrt{3g_2}, \quad \lambda_3=3e_2,\quad \lambda_2=-\lambda_4,\quad \lambda_1=3e_3.
$$
Here $g_2=4(e_2^2+e_3^2+ e_2e_3)$. The shift by $\lambda_5$ brings the genus two hyperelliptic curve with branch points $\lambda_1, \dots, \lambda_5$ to the curve of the form \eqref{hyperx} with  branch points at $0$, $x_1=3e_2+6e_3$, $x_2=6e_2+3e_3$, $u_1= \sqrt{12(e_2^2+e_3^2+ e_2e_3)}+ 3(e_2+e_3)$ and $u_2= -\sqrt{12(e_2^2+e_3^2+ e_2e_3)}+ 3(e_2+e_3)$. It is immediate that
$$
e_3=\frac{2x_1-x_2}{9},\qquad e_2=\frac{2x_2-x_1}{9}.
$$
This example is setting stage for the following Proposition.
\end{example}

\begin{proposition}\label{prop:2gaplam} Given the curve $\Gamma_{e_2,e_3}$ \eqref{GW} with $e_1< e_2 < e_3$, denote $x=e_2+2e_3$ and $u=2e_2+e_3$. There exists a unique continuous deformation of the curve $\Gamma_{e_2,e_3}$ \eqref{GW} which preserves the real period $2w_1$ of the potential $v_2$, such that  $u=u(x)$ is the real-valued function of $x$.   Under this deformation, the genus two hyperelliptic curve $\mathcal C_{x_1,x_2}$ of the form \eqref{hyperx} with the branch points at $0$, $x_1=3e_2+6e_3$, $x_2=6e_2+3e_3$, $u_1= \sqrt{12(e_2^2+e_3^2+ e_2e_3)}+ 3(e_2+e_3)$ and $u_2= -\sqrt{12(e_2^2+e_3^2+ e_2e_3)}+ 3(e_2+e_3)$ deforms in a way that  $x_1=3x$,  $x_2=3u(x_1/3)$  and $u_1=u_1(x_1, x_2)$ and  $u_2=u_2(x_1, x_2)$ are real-valued functions  of real variables that satisfy equations from Example \ref{examplegenus2} and \eqref{derivative}  with $\Omega$ and $v_k$  given by \eqref{Omega} and \eqref{v} on the surface $C_{x_1,x_2}$.
\end{proposition}

{\it Proof.}
According to Corollary \ref{prop:onegap}, there exists a unique continuous deformation of the curve $\Gamma_{e_2,e_3}$ \eqref{GW} which preserves the period $2w_1$ of the potential $v_1(X)=2\wp(X)$, such that $u=2e_2+e_3$ satisfies equation \eqref{ode} as a function of $x=2e_3+e_2$ with the initial condition \eqref{uprime},  where $\Omega=\Omega_0$ on $\Gamma_{e_2,e_3}$ is from the proof of Theorem \ref{thm_2w1}, see \eqref{Omega0}. Due to  Corollary \ref{cor_real}, this deformation has a unique restriction to a real neighborhood $\mathcal X$ of $x$, for which $u=u(x)$ is a real-valued function. Using  Example \ref{ex:2gap}, this generates a deformation  of the  genus two hyperelliptic curve  $\mathcal C_{x_1,x_2}$ \eqref{hyperx} with  $x_1, \;x_2, \; u_1,\; u_2$ as in the statement of the proposition, where  $x_1=3x$ and $x_2=3u(x_1/3)$ are both real   and where the associated periodic potential is $v_2$.  Note now that, due to the fact that $e_1, e_2, e_3$ are real,  the $\wp$-function has  a real  and a purely imaginary  periods, denoted by $w_1$ and $w_2$ respectively (see e.g. \cite{WW}, Example 1, Ch.20-32, p. 444). Thus, the $\wp$-function  is real-valued along the boundary of the fundamental rectangle with the vertices $0$, $w_1$, $w_1+w_2$, $w_2$ (see e.g. \cite{WW}, Example 2, Ch.20-32, p. 444). Thus, given that the curve $\mathcal C_{x_1,x_2}$ for $x\in\mathcal X$ is associated to a real-periodic potential, it is a real Hill curve \cite{Mei, MM}, see also the proof of Corollary \ref{cor_real} and thus the deformation in question is an isoperiodic deformation of the pair $(\mathcal C_{x_1,x_2}, \Omega_0)$, where $\Omega$ is defined in \eqref{Omega} for $\alpha=0$ on  $\mathcal C_{x_1,x_2}$ and satisfies the Hill condition \eqref{eq:hill} with $T=2w_1$. By Theorem \ref{thm_main}, the functions $u_1=u_1(x_1, x_2)$ and  $u_2=u_2(x_1, x_2)$  satisfy equations from Example \ref{examplegenus2} and \eqref{derivative}, with  $\Omega=\Omega_0$. The reality of $u_1$ and $u_2$ through the deformation follows from reality of the endpoints of the gaps for a real-periodic potential.
$\Box$\\

\section{Isoperiodic deformations of solutions to the Neumann problem}\label{sec:neum}

 In this section we apply Theorems \ref{thm_main}, \ref{thm_converse-g} and Corollary \ref{cor_real} to obtain isoperiodic deformations  of periodic solutions  to the general Neumann system on $S^n$.
The Neumann system from classical mechanics describes motion of a particle on an $n$-dimensional sphere
$$
S^n:\quad \sum_{j=1}^{n+1} q_j^2=1
$$
subject to the potential force with the quadratic potential
$$
V(q)=\frac{1}{2}\sum_{j=1}^{n+1}A_jq_j^2, \quad A_j = \text{const}.
$$
The equations of motion are
$$
\ddot q_j=-A_jq_j+F(t)q_j, \quad j=1, \dots, n+1.
$$
Here $F(t)=\sum_{j=1}^{n+1}(A_jq_j^2-\dot q_j^2)$ is the Lagrange multiplier corresponding to the reaction of the constraint which forces the particle to stay on the sphere $S^n$. This system is Hamiltonian restricted to the sphere $S^n$ with the Hamiltonian
$$
H=\frac{1}{2}\Big(\sum_{j=1}^{n+1}A_jq_j^2+q^2p^2-(qp)^2\Big),
$$
where $q=(q_1, \dots, q_{n+1}), \;p=(p_1, \dots, p_{n+1})$ and $qp$ is the euclidean scalar product.
The independent involutive first integrals of motion are
$$
F_k(q, p)=q_k^2+\sum_{j\ne k}\frac{(q_kp_j-q_jp_k)^2}{A_j-A_k}, \quad k=1, \dots, n+1.
$$
Thus, the Neumann system is a completely integrable system in the Liouville sense and $2H=\sum_{j=1}^{n+1}A_jF_j$. The case $n=1$ was considered in Section \ref{sec:onegap}, see Example \ref{ex:planarneum}.
For  positive $A_j$ and $n\geqslant 2$, the Neumann system is related to another famous integrable system, the Jacobi problem of geodesics on the ellipsoid in $\mathbb R^{n+1}$
$$
\sum_{j=1}^{n+1}\frac{\tilde q_j^2}{A_j}=1.
$$
The correspondence between the Neumann and the Jacobi problems is as follows
$$
\tilde q=p, \quad \tilde p=-q, \quad\tilde H=\sum_{j=1}^nA_j^{-1}F_j.
$$

There is a remarkable connection between the Neumann system on $S^n$ and $n$-gap potentials, according to \cite{Mo}, \cite{Mo1}, \cite{Ve}, \cite{Mu}. Consider an $n$-gap differential operator
\begin{equation}
\label{LL}
\mathbb L=-\frac{d^2}{dx^2}+v(x),
\end{equation}

where the potential $v(x)$ is periodic. Let the curve $\Gamma_n$ given by the equation $w^2=P_{2n+1}(z)$ be the spectral curve of $\mathbb L$, that is, denoting $z_{2n+1}<\dots<z_1$ the zeros of the polynomial $P_{2n+1}$, we have that the Bloch spectrum of $\mathbb L$ is composed of the intervals $[z_{2j+1}, z_{2j}]$ for $j=1, \dots, n$ and $[z_1, \infty].$ Denote by $P_{z_k}$ the point of  $\Gamma_n$ for which $z=z_k$.
Assume that the $a$- and $b$-cycles of the spectral curve $\Gamma_n$ are selected as in Example \ref{ex:realhyp}.

 Denote also $-A_j=z_{2j+1}$ for $j=1, \dots, n$. The eigenvectors of $\mathbb L$ corresponding to eigenvalues $-A_j$ are constructed by considering  the Baker-Akhiezer function $\psi(x, P)$ with $P\in\Gamma_n$ and $x\in\mathbb R$ corresponding to the operator $\mathbb L.$ Setting
$$
\psi_j(x)=\kappa_j\psi(x, P_{z_{2j+1}}),
$$
with
$$
\kappa_j=\Big(\prod_{k\ne j}(A_j-A_k)\Big)^{-\frac{1}{2}}, \quad j=1, \dots, n+1,
$$
we have
$$
\mathbb L\psi_j=-A_j\psi_j\,.
$$
This last equation is equivalent to the Neumann system
$$
\ddot \psi_j=-A_j\psi_j + v(x)\psi_j,\quad j=1, \dots, n+1,
$$
with the following correspondence $x\rightarrow t$, $\psi_j\rightarrow q_j$, $v(x)\rightarrow F(t)$, taking into account that, see \cite{Ve},
$$
\sum_{j=1}^{n+1}\psi_j^2=1.
$$
The original classical Neumann system is for $n=2$ and corresponds to two-gap potentials.
The solutions of the classical Neumann system are given by the formulae:
\begin{align}
\label{qsol}
\nonumber
q_1(t)&=\kappa_1\frac{\theta[(0, 1/2),(0,0)](tU+z_0)\theta(z_0)}{\theta[(0, 1/2),(0,0)](z_0)\theta(tU+z_0)},\\
q_2(t)&=\kappa_2\frac{\theta[(1/2, 0),(0, 1/2)](tU+z_0)\theta(z_0)}{\theta[(1/2, 0),(0, 1/2)](z_0)\theta(tU+z_0)},\\
\nonumber
q_3(t)&=\kappa_3\frac{\theta[(0,0),(1/2, 1/2)](tU+z_0)\theta(z_0)}{\theta[(0,0),(1/2, 1/2)](z_0)\theta(tU+z_0)}.
\end{align}
 Here $\theta$ is the theta function associated with the real  curve $\Gamma_2$, and
 $\kappa_j$ are defined above and $U$ \eqref{U} is the same as in Section \ref{sec:twogap}  for the curve $\Gamma_2$.
The solutions of a general Neumann system on $S^n$ for arbitrary $n$ have practically the same form in terms of $\theta$-functions related to genus $n$ hyperelliptic real  curves $\Gamma_{n}$, see \cite{Mu}.
Solutions of the Neumann system are periodic with period $T\in \mathbb R$ if and only if the corresponding $n$-gap potential $v(x)$ is periodic, thus if and only if the period vector $U$ of $\Omega_0$ satisfies \eqref{eq:periodicity} with some $T\in \mathbb R$. This is  true if and only if the Hill  condition \eqref{eq:hill} is satisfied. Note that all $m_j$ are zero in \eqref{eq:periodicity}, since $U$ is a real vector. The reality of $U$ follows from the fact that the relevant hyperelliptic curves are real, see Example \ref{ex:realhyp}, the reality of the coefficients of the normalized holomorphic differentials and formula \eqref{Omega0g} expressing $U$ in terms of those differentials.

In the case $n=2$, by a rescaling of time $t$ we may assume that $A_3=0$. (Similarly, we may assume that $A_{n+1}=0$ for general $n$.)
Proposition \ref{prop:2gap}, implies the following statement.

\begin{proposition}
\label{prop:Neum}
Let $q=(q_1,q_2,q_3)$ be a periodic  solution \eqref{qsol} to the classical Neumann system for $n=2$ with the potential $V(q)= \frac{1}{2}(A_1q_1^2+A_2q_2^2)$ and $A_1,\,A_2\in\mathbb R.$ Denote the period of the solution by $T\in \mathbb R$. Let $\Gamma$ be the hyperelliptic curve of genus $2$ represented as a two-fold ramified covering of $z$-sphere with a branch point at $z=\infty$ and finite branch points $z_5=-A_3=0<z_4<z_3=-A_2<z_2<z_1=-A_1$, the spectral curve of the system. Assume now that $\{\Gamma_A, q_A\}$ is a continuous family of hyperelliptic curves with finite branch points $z_5=-A_3=0<z_4<z_3=-A_2<z_2<z_1=-A_1$ and denote  the corresponding solutions by $q_A$ \eqref{qsol}.   Denote $x_1=-A_1$, $\;x_2=-A_2$, $\;u_1=z_2$, and $u_2=z_4$. Let  $x_1$ and $x_2$ (that is $A_1$ and $A_2$) vary in some real neighborhood. The obtained solutions $q_A$ are periodic with the same period $T$ if and only if $u_1, u_2$ as real-valued functions of  $x_1, x_2$ satisfy equations of Example \ref{examplegenus2} and the initial conditions \eqref{derivative}. Here $\Omega$ and $v_j$ in  \eqref{derivative} are differentials on $\Gamma_A$ defined by \eqref{Omega} and \eqref{v} with $\lambda=z$, respectively.
 More generally, assume that $\{\Gamma_A, q_A\}$ is a continuous family of hyperelliptic curves $\Gamma_A$ of genus $n$ represented as a two-fold ramified coverings of the $z$-sphere with finite branch points $z_{2n+1}<\dots<z_1$, with $z_{2j-1}=-A_j$, $j=1, \dots, n+1$ where $A_{n+1}=0$. Denote by $q_A$ the  solution  to the general Neumann system on $S^n$ with the potential $V(q)= \frac{1}{2}\sum_{j=1}^{n}A_jq_j^2$ constructed  as in \cite{Mu} from the corresponding spectral curve $\Gamma_A$, which is periodic with a period $T\in \mathbb R$. Denote $x_j=z_{2j-1}=-A_j$, $u_j=z_{2j}$, $j=1,\dots, n$. By varying $A_1, \dots, A_n$ in some real neighborhoods, the obtained solutions $q_A$ are periodic with the same period $T$ if and only if
 $u_1,\dots, u_n$ as real-valued functions of  $x_1,\dots, x_n$ satisfy equations \eqref{umxnxk} and \eqref{umxkxk} of Theorem \ref{thm_main} and the initial conditions \eqref{derivative}. As before, $\Omega$ and $v_j$ in  \eqref{derivative} are differentials on $\Gamma_A$ defined by \eqref{Omega} and \eqref{v} with $\lambda=z$, respectively.
\end{proposition}

{\it Proof.} A periodic solution to the Neumann system corresponds to a periodic finite-gap potential \cite{Mo,Mo1,Ve}, and thus the corresponding spectral curve   has real branch points  and satisfies the  Hill condition \eqref{eq:hill} \cite{MM}.   Therefore, the curves $\Gamma_A$ in the family are  curves with real branch points, satisfying the  Hill condition \eqref{eq:hill}. The statement thus follows from Theorem \ref{thm_main}, Theorem \ref{thm_converse-g} and Corollary \ref{cor_real}.
$\Box$

\section{Isoperiodic deformations of the finite-gap periodic solutions to the KdV equation}
\label{sect_kdv}

 In this section, we obtain isoperiodic deformation of periodic solutions to the KdV equation as a consequence of  Theorems \ref{thm_main} and \ref{thm_converse-g} and Corollary \ref{cor_real} from Section \ref{sect_main}.

 We consider a family $\mathcal C_\x$ of compact genus $g$ hyperelliptic Riemann surfaces  corresponding to the curves of equation \eqref{hyperx}. The hyperelliptic coverings of the $\lambda$-sphere have a ramification point $P_{\infty}$ over $\lambda=\infty.$
 Denoting, as before, the other ramification points of the hyperelliptic covering by $P_0,$ $\puj$ and $\pxj$,
%let us assume there are branch cuts between $P_0$ and $P_\infty$ as well as between $\pxj$ and $\puj$ for $j=1, \dots, g.$
let $\{\a_1,\dots, \a_g; \b_1, \dots, \b_g\}$ be a fixed canonical homology basis on the surfaces $\mathcal C_\x$ for $\x \in \Sx$, such that $\a_j$ goes around the points $\pxj$ and $\puj$ while  $\b_j$ goes around the points $P_0, \pxk, \puk$ and $\pxj$ for $k=1, \dots, j-1$. The cycles are oriented in a way so that the intersection indices be $\a_k\circ \a_j=\b_k\circ \b_j=0$ and $\a_k\circ \b_j=\delta_{kj}$. The domain $\Sx$ is chosen as in Section \ref{sect_iso}.

Consider  differentials of the second kind $\Omega^{(n+1)}_{P_{\infty}}(P)$ with $P\in \mathcal C_\x$ normalized by  vanishing of their $a$-periods and having a single pole at $P_{\infty}$ of order $n+2$ with the local behaviour in terms of the local parameter $\zeta_\infty$ \eqref{coordinates}
\begin{equation}
\label{eq:normsecond}
\Omega^{(n+1)}_{P_{\infty}}(P)  = \left( \frac{1}{\zeta^{n+2}_\infty(P)} + \mathcal O(1) \right)  \,,\quad P\sim P_\infty,
\end{equation}

This way, for the differential $\Omega_0$ defined by \eqref{Omega} with $\alpha=0$, we have $\Omega_0=\Omega^{(1)}_{P_{\infty}}$. Denote further
\begin{equation}\label{eq:UV}
U=-\frac{1}{2\pi\i}\left(\oint_{b_1}\Omega^{(1)}_{P_{\infty}}, \dots, \oint_{b_g}\Omega^{(1)}_{P_{\infty}}\right)^T\quad\text{and}\quad
W= - \frac{3}{2\pi\i}\left(\oint_{b_1}\Omega^{(3)}_{P_{\infty}}, \dots, \oint_{b_g}\Omega^{(3)}_{P_{\infty}}\right)^T.
\end{equation}
We follow  \cite{Dub1981} closely in this section while adjusting some constants to continue using our normalization \eqref{normalization}, \eqref{bW}.
For a fixed $\x$, the meromorphic function $\lambda:\mathcal C_\x\to\mathbb CP^1$  has a unique pole which is of the second order at $P_{\infty}$  on $\mathcal C_\x$. With the help of this function, one obtains a reduction of the Kadomtsev-Petviashvili equation (KP) related to the Baker-Akhiezer function $\phi$ with an essential singularity at $\P_\infty$ defined by the polynomial
$$
q(k)=kX+k^3t,
$$
where $k$ is the reciprocal of the local parameter $\zeta_\infty$ around $P_{\infty}$ on $\mathcal C_\x$.
%Thus, we may assume that $\lambda=k^2$.
The Baker-Akhiezer function at $P\sim P_\infty$ behaves as
$$
\phi(X, t;P)=e^{kX+k^3t}\big(1+\frac{\xi_1}{k}+\frac{\xi_2}{k^2}+\dots \big),
$$
where $\xi_i$, $i=1,2$ are functions of $X$ and $t$ to be determined below.
Consider the operators
$$
\mathbb L=\partial^2_{X}+v, \quad \mathbb A=\partial^3_{X}+\frac{3}{2}v\partial_X + w,
$$
where
$$
v=-2\frac{\partial \xi_1}{\partial X},\quad w=3\xi_1\frac{\partial \xi_1}{\partial X}+3\frac{\partial^2 \xi_1}{\partial X^2}-3\frac{\partial \xi_2}{\partial X}.
$$
Then the Baker-Akhiezer function $\phi$ satisfies the system
$$
\mathbb L\phi=\lambda \phi, \qquad \frac{\partial \phi}{\partial t}=\mathbb A\phi.
$$
The compatibility condition for this system leads to the L-A equation
$$
\frac{\partial}{\partial t}\mathbb L=\big[\mathbb L, \mathbb A\big],
$$
which is equivalent to the system of nonlinear PDEs for $v$ and $w$:
$$
\frac{3}{4}v_{XX}-w_X=0,\qquad -v_t+v_{XXX} +\frac{3}{2}vv_X-w_{XX}=0.
$$
By excluding $w$ from the last two equations, one gets the following form of the Korteweg-de Vries (KdV) equation \cite{Dub1981}:
\begin{equation}\label{eq:kdv}
v_t=\frac{1}{4}(6vv_X+v_{XXX}).
\end{equation}
 As indicated in e.g. \cite{Dub1981}, the function
\begin{equation}\label{eq:solkdv}
v(X, t)=2\frac{\partial^2}{\partial X^2}\ln \theta(XU+tW+z_0)+c,
\end{equation}
solves the KdV equation \eqref{eq:kdv}, where $\theta(z)$ is the theta-function associated with $\mathcal C_\x$ and $U$ and $W$ are the $b$-periods of the differentials of the second kind \eqref{eq:UV} \eqref{eq:normsecond} with a unique pole at $P_{\infty}$. Considering now a variation of $\x\in\Sx$, as an application of Theorems \ref{thm_main} and \ref{thm_converse-g}, we have the following statement.
\begin{proposition}
\label{thm:kdvappl}
Let $\mathcal C_{\x_0}$ be a compact Riemann surface of the curve \eqref{hyperx} with $\x=\x_0$ and let $v=v_{\x_0}$ be a solution \eqref{eq:solkdv} of the KdV equation \eqref{eq:kdv}, constructed from $\mathcal C_{\x_0}$ and the point $P_{\infty}\in \mathcal C_{\x_0}.$ Then there exists a unique continuous family of surfaces $\mathcal C_\x$ parametrized by $\x\in\mathcal X$ with $\x_0\in\mathcal X$ such that all solutions from the corresponding continuous family of solutions $v_{\x}$ of the KdV equation  obtained from $\mathcal C_{\x}$ and $P_{\infty}\in \mathcal C_{\x}$  by \eqref{eq:solkdv} have the same   vector  $U$  (the wavevector). Moreover,  for the family $\mathcal C_{\x}$, the functions $u_1(\x), \dots, u_g(\x)$ from \eqref{hyperx}  satisfy equations \eqref{umxnxk}, \eqref{umxkxk} of Theorem \ref{thm_main}  and  \eqref{derivative}, where $\Omega$ in \eqref{derivative} is $\Omega_0=\Omega^{(1)}_{P_{\infty}}$ and $v_j$ are defined by \eqref{v}.
\end{proposition}
{\it Proof.} The family $\mathcal C_\x$ is an isoperiodic family $(\mathcal C_\x, \Omega_{P_\infty}^{(1)})$ relative to the normalized differential $\Omega_{P_\infty}^{(1)}$ \eqref{eq:normsecond}.  The existence and uniqueness of such a family is given by Theorem \ref{thm_converse-g}. The conditions \eqref{derivative}  come from Theorem \ref{thm_derivatives}. The statement concerning equations of Theorem \ref{thm_main} holds due to the fact that these equations  are derived starting from the condition that the $b$-periods of the differential $\Omega$ \eqref{Omega} are constant for  all surfaces from the family $\mathcal C_\x$. The statement of Theorem \ref{thm_main} applies to the normalized differentials $\Omega_{\pinfty}^{(1)}$ and their $b$-periods  equal to $-2\pi\i U$ \eqref{U}, by setting $\alpha=0$ in \eqref{Omega}.
$\Box$

\bigskip
If there exists a real nonzero number $T$, such that $TU$  belongs to the lattice of the Jacobian of $\mathcal C_\x$ that is if \eqref{eq:periodicity} is satisfied for some integers $n_j$ and $m_j$, then $v$ \eqref{eq:solkdv}  is periodic in $X$ with a period $T$.  In this case, we say that $v$ is a {\it spatial $T$-periodic function}.
%The lattice of the Jacobian of $\mathcal C_\x$  is generated by $(Id, \mathbb B)$, where $Id$ is the identity matrix and $\mathbb B$ is the Riemann matrix, see \eqref{normalization} and \eqref{B}.

It follows from  \cite{DMN},  that all real smooth finite-gap solutions of the KdV equation are obtained by \eqref{eq:solkdv} from a hyperelliptic curve $\Gamma$ \eqref{eq:gammareal} of genus $g$ which has all branch points real $\lambda_{2g+1}<\lambda_{2g}<\dots< \lambda_1$,  with the canonical basis of cycles chosen as in Exemple \ref{ex:realhyp}, and with an additional condition as follows. It is required that the non-special degree $g$ divisor $\mathcal D=\sum_{j=1}^g(\gamma_j, w_j)$ of poles of the Baker-Akhiezer function $\phi$ be such that
$\lambda_{2j}\le \gamma_j\le \lambda_{2j-1}$, for $j=1,\dots g$. The divisor $\mathcal D$  determines the constant $z_0$ in \eqref{eq:solkdv} up to the vector of Riemann constants $K$ of $\Gamma$ (note that $\mathcal D=Q_1+Q_2$ in \eqref{eq:z0} for $g=2$). In this case, $U$ is real and thus such a solution $v$ \eqref{eq:solkdv} is spatially periodic with the period $T$ if and only if the curve $\Gamma$ satisfies the  Hill condition  \eqref{eq:hill}.

\begin{proposition}
\label{thm:kdv}
Given  $v_{\x_0}=v_{\x_0}(X,t)$, an  arbitrary smooth real finite gap solution  \eqref{eq:solkdv} of the KdV equation \eqref{eq:kdv} which is spatially periodic with the period $T$, denote by  $\mathcal C_{\x_0}$ a curve  \eqref{hyperx} with $\x=\x_0$ and the point $P_{\infty}\in \mathcal C_{\x_0}$, from which   $v_{\x_0}$ is constructed. There exists a continuous family of  solutions $v_{\x}(X, t)$ for $\x\in \mathcal X\cap \mathbb R^g$ with $\mathcal X$ being some open neighbourhood of $\x_0$ for which all solutions $v_{\x}(X,t)$ of the KdV equation  are spatially periodic with period $T$. Moreover, there exists a family of curves $\mathcal C_{\x}$, $\x\in \mathcal X\cap \mathbb R^g$, given by \eqref{hyperx}, which satisfy   the Hill condition  \eqref{eq:hill}  relative to a canonical homology basis chosen as in Example \ref{ex:realhyp}. The family of curves $\mathcal C_{\x}$ is  such that $v_{\x}(X,t)$ is constructed by \eqref{eq:solkdv} from $\mathcal C_\x$ and $P_{\infty}\in \mathcal C_\x$ and the real-valued functions $u_1(\x), \dots, u_g(\x)$ from \eqref{hyperx} satisfy equations \eqref{umxnxk}, \eqref{umxkxk} of Theorem \ref{thm_main} as well as \eqref{derivative}.  The differentials in \eqref{derivative} are $\Omega_0=\Omega^{(1)}_{P_{\infty}}$ \eqref{eq:normsecond} and $v_j$  \eqref{v} defined on the curves $\mathcal C_\x$;  $\Omega_0$ is normalized by vanishing of its periods with respect to the $a$-cycles invariant under the antiholomorphic involution.
\end{proposition}
{\it Proof.} From the fact that  $v_{\x_0}$  is a real smooth solution of KdV equation, it follows  that the covering $\lambda:\mathcal C_{\x_0}\to\mathbb CP^1$ has all branch points real. Therefore, we may choose a canonical homology basis on $\mathcal C_{\x_0}$ to satisfy relations $\rho a_j=a_j$, $\rho b_j=-b_j$ with respect to the antiholomorphic involution $\rho$ as in Example \ref{ex:realhyp}. From the periodicity of $v_{\x_0}$ and from the properties of quasi-periodicity of the theta-function, it also follows that the curve $\mathcal C_{\x_0}$ satisfies the  periodicity   condition \eqref{eq:periodicity}  with  $T$ and $m_j=0$, $j=1, \dots, g.$   The integers $m_j$ in \eqref{eq:periodicity} vanish due to the definition \eqref{eq:UV} of the vector $U$ and its reality; recall that the Riemann matrix is purely imaginary for this choice of homology basis. Theorem \ref{thm_converse-g} gives the existence and the uniqueness of an isoperiodic family $(\mathcal C_\x, \Omega^{(1)}_{P_\infty}(\x))$ where $\Omega^{(1)}_{P_\infty}(\x)$ is the normalized differential of the second kind \eqref{eq:normsecond} on $\mathcal C_\x$ and where $\x\in\mathcal X$ for some open neighbourhood $\mathcal X$ of $\x_0$. Given that this family is isoperiodic, after restricting it to $\x\in\mathcal X\cap \mathbb R^g$, we obtain that the branch points of the obtained coverings  $\lambda:\mathcal C_\x\to\mathbb CP^1$ are all real (the functions  $u_j(\x)$ are real due to Corollary \ref{cor_real}). For this restriction, the homology basis chosen on $\mathcal C_{\x_0}$ extends to a basis on $\mathcal C_{\x}$ satisfying the same conditions with respect to the anitholomorphic involution on $\mathcal C_{\x}$ with $\x\in\mathcal X\cap \mathbb R^g$, that is $\rho a_j=a_j$, $\rho b_j=-b_j$. The differential $\Omega_0=\Omega^{(1)}_{P_\infty}(\x)$ is thus normalized using this basis and  the   Hill  condition \eqref{eq:hill} holds  for $\Omega_0=\Omega^{(1)}_{P_\infty}(\x)$ with the same $T$ for all $\x\in \mathcal X\cap \mathbb R^g$  and the above homology basis. Thus \eqref{eq:periodicity} is satisfied with the same $T$ and $m_j=0$, $j=1, \dots, g,$ for all  the curves $\mathcal C_\x, \; \x\in\mathcal X\cap \mathbb R^g.$ As a consequence, solutions $v_\x(X,t)$ constructed from $\mathcal C_\x$ are spatial periodic with period $T.$ This proves the existence of a continuous family of solutions to the KdV sharing the same spatial period $T$. Equations  \eqref{umxnxk}, \eqref{umxkxk} of Theorem \ref{thm_main} are satisfied since the family $(\mathcal C_\x, \Omega^{(1)}_{P_\infty}(\x))$ is isoperiodic, that is since the $b$-period vector $U$ is preserved over the family.  The derivatives are expressed by  \eqref{derivative}, due to Theorem \ref{thm_derivatives} .
$\Box$

\section{Isoperiodic deformations and the sine-Gordon equation}
\label{sect_sG}

 In this section we discuss one more application of Theorems \ref{thm_main} and \ref{thm_converse-g} and Corollary \ref{cor_real} from Section \ref{sect_main}, an application to the theory of periodic solutions of the sine-Gordon equation.
Here, the family of curves $\mathcal C_\x$ for $\x \in \Sx$, their canonical homology bases and ramification points of the coverings $\lambda:\mathcal C_\x\to\mathbb CP^1$ are as in Section \ref{sect_kdv}.

Following \cite{DubNat}, let us pick two ramification points of the covering $\lambda$ and denote them by $\hat P$ and $\hat Q$. Let $p$ and $q$  be local parameters at $\hat P$ and $\hat Q$, respectively, proportional to the respective standard local parameters \eqref{coordinates} $\zeta_{\lambda_k}$, where the proportionality coefficients satisfy a precise condition from \cite{DubNat} unimportant for the present discussion. Denote by $U=(U_1, \dots, U_g)$ and  $V=(V_1, \dots, V_g)$ the vectors
defined by
\begin{equation}\label{eq:UVsg}
U_j=\omega_j(\hat P), \quad V_j=\omega_j(\hat Q).
\end{equation}
Here $\omega_1, \dots, \omega_g$ is a basis of normalized holomorphic differentials, as before,  and the evaluation at the points $\hat P$ and $\hat Q$ is done with respect to the local parameters $p$ and $q$, respectively, as in \eqref{evaluation}.
Thus,  $2\pi\i\,U$ and $2\pi\i\,V$ are the vectors of $b$-periods of differentials $\Omega_{\hat P}$ and $\Omega_{\hat Q}$ of the second kind with poles of the second order only at $\hat P$ and only at $\hat Q$, respectively, and with the leading coefficients of the Laurent series in $p$ and $q$ being equal to one.
Define
$$
\Delta=\big(\int_{\hat Q}^{\hat P}\omega_1,\dots,  \int_{\hat Q}^{\hat P}\omega_g  \big)^T.
$$
Then, since $\hat P$ and $\hat Q$ are ramification points, $\Delta$ is a vector of  half-periods:
$$
\Delta=\frac{1}{2}M+\frac{1}{2}\mathbb BN, \quad M, N\in \mathbb Z^g,
$$
where $\mathbb B$ is the Riemann  matrix. According to \cite{DubNat}, the function
\begin{equation}
\label{eq:solSG}
u(X, Y)=\frac{1}{\i}\Big[\ln\frac{\theta(XU+YV+\Delta+z_0)\theta(XU+YV-\Delta+z_0)}{\theta^2(XU+YV+z_0)} +\pi\i\langle N, \Delta\rangle\Big],
\end{equation}
is a solution of the sine-Gordon equation
\begin{equation}\label{eq:SG}
u_{XY}=-4\kappa \sin u,
\end{equation}
where
$$
\kappa = \exp (-\pi\i\langle N, \Delta\rangle),
$$
and  $z_0$ is an arbitrary vector in the Jacobian of $\mathcal C_\x$.
According to \cite{DubNat},  real solutions to the sine-Gordon equation correspond to real hyperelliptic curves,  that is  curves defined by real polynomials  and the ramification points $\hat P, \hat Q$ belong to the same real oval,  as defined in Example \ref{ex:realhyp}. In this case $\Delta$ is real,
$$
\Delta=\frac{1}{2}M, \quad M\in\mathbb Z^g,
$$
and solution \eqref{eq:solSG} simplifies to
\begin{equation}
\label{eq:solSGreal}
u(X, Y)=-\i\ln\left[\frac{\theta(XU+YV-\frac{M}{2}+z_0)}{\theta (XU+YV+z_0)}\right]^2,
\end{equation}
while the sine-Gordon equation simplifies to
\begin{equation}
\label{eq:SG1}
u_{XY}=-4\sin u.
\end{equation}
 Solution $u(X, Y)$ given by \eqref{eq:solSGreal} of equation  \eqref{eq:SG1} is real \cite{DubNat}  under some condition on the vector $z_0$ (see two paragraphs below).

 Our isoperiodicity results in this paper (Theorems \ref{thm_main}, \ref{thm_converse-g},  Corollary \ref{cor_real}) are formulated for  a differential $\Omega$ \eqref{Omega} of the second kind having a pole of the second order at $P_\infty$, the ramification point at infinity of $\mathcal C_\x$. These results can be reformulated in a straightforward, though tedious way for  a pole at any other ramification point. For simplicity  and without loss  of generality, we will assume from now on, that $\hat P=P_{\infty}$.

Let us assume now that we have a real periodic solution $\hat u_{\x_0}$ constructed from a real hyperelliptic curve  \eqref{hyperx}, where {\it all
the branch points $x^0_j, u^0_j$, $j=1, \dots g$ are real} (note that not all real finite-gap solutions of the sine-Gordon equation satisfy that).
 The homology basis chosen at the beginning of the section, see Section \ref{sect_kdv}, satisfy relations  from Example \ref{ex:realhyp} with respect to the antiholomorphic involution on this curve.
 Here we have that  $\hat Q$ corresponds to the largest real finite  branch point, so that it belongs to the same real oval as $P_{\infty}$, see Example \ref{ex:realhyp}. In this case, for the reality of solutions, it is also required that $\bar z_0\equiv -\Delta-z_0$.  (Here the sign $\equiv$ means modulo the Jacobian lattice.)

 Assume that $x^0_g$ is strictly larger than all other finite branch points $u^0_1, \dots u^0_g, x^0_1, \dots, x^0_{g-1}, 0$, and that $\hat Q=P_{x^0_g}$ is the ramification point corresponding to $x^0_g$.  Since the  canonical basis of cycles satisfies conditions of Example \ref{ex:realhyp}, the corresponding holomorphic  normalized differentials   $\omega_1, \dots, \omega_g$ are all real. The vector $U$ from \eqref{eq:UVsg} is also real, and the solution   $\hat u_{\x_0}$ given by  \eqref{eq:solSGreal} to the sine-Gordon equation \eqref{eq:SG1} is  periodic in $X$ if and only if there exists a positive $T$ such that $TU$ is an integer-valued vector, $TU\in \mathbb Z^g$.
\begin{proposition}
\label{thm_SGperiodic}
Let $\mathcal C_{\x_0}$ be a  curve of equation \eqref{hyperx} with $\x=\x_0$  with all branch points of the covering $\lambda$ being real. Let $\hat u_{\x_0}=\hat u_{\x_0}(X,Y)$ be the real  solution  \eqref{eq:solSGreal} to the sine-Gordon equation \eqref{eq:SG1} constructed from the surface $\mathcal C_{\x_0}$   and from the ramification points $\hat P=P_{\infty}\in \mathcal C_{\x_0}$ and $\hat Q$  such that $\lambda(\hat Q)=\max\{0, x^0_j, u^0_j|, j=1,\dots, g\}$. Assume that $\hat u_{\x_0}$ is periodic in $X$ with the period $T$.  There exists a continuous family of real solutions $\hat u_{\x}(X, Y)$ for $\x\in \mathcal X\cap \mathbb R^g$ with $\mathcal X$ being some open neighbourhood of $\x_0$ for which all solutions $\hat u_{\x}(X, Y)$ of the sine-Gordon equation  \eqref{eq:SG1} are  periodic with the period $T$. Moreover, there exists a family of curves $\mathcal C_{\x}$, $\x\in \mathcal X\cap\mathbb R^g$, given by \eqref{hyperx} such that $\hat u_{\x}(X, Y)$ is   constructed by \eqref{eq:solSGreal} from $\mathcal C_\x$, $\hat P=P_{\infty}\in \mathcal C_\x$ and $\hat Q$ with   $\lambda(\hat Q)=\max\{0, x_j, u_j|, j=1,\dots, g\}$. The real-valued functions $u_1(\x), \dots, u_g(\x)$  from \eqref{hyperx} corresponding to the family $\mathcal C_\x$, $\x\in \mathcal X\cap \mathbb R^g$, satisfy equations \eqref{umxnxk}, \eqref{umxkxk} of Theorem \ref{thm_main}.
\end{proposition}
\medskip
{\it Proof.}  We have that the curve $\mathcal C_{\x_0}$ produces a real periodic solution  $\hat u_{\x_0}=\hat u_{\x_0}(X,Y)$ \eqref{eq:solSGreal} of the sine-Gordon equation and that all branch points of the covering $\lambda:\mathcal C_{\x_0}\to\mathbb CP^1$ are  real.  In particular, the canonical homology basis satisfies $\rho a_j=a_j$, $\rho b_j=-b_j$ where $\rho$ is the antiholomorphic involution on $\mathcal C_{\x_0}$. Thus the vector $U$  \eqref{eq:UVsg} is real and the Riemann matrix $\mathbb B$ is imaginary. From this and  the periodicity of the solution $\hat u_{\x_0}$ \eqref{eq:solSGreal} and from the properties of quasi-periodicity of the theta-function, we have that the vector $U$ satisfies \eqref{eq:periodicity} with $m_j=0$ and $T>0$ being the period of $\hat u_{\x_0}$. This implies that the
curve $\mathcal C_{\x_0}$   satisfies the  Hill condition \eqref{eq:hill}, where $\Omega_0=\Omega^{(1)}_{P_\infty}$ is the differential \eqref{eq:normsecond} normalized with respect to the above homology basis.  Theorem \ref{thm_converse-g} and Corollary \ref{cor_real} give the existence of an isoperiodic real family $(\mathcal C_\x, \Omega^{(1)}_{P_\infty})$ for $\x\in\mathcal X\cap \mathbb R^g$ for some open neighbourhood $\mathcal X$ of $\x_0$ providing isoperiodic deformations of the pair $(\mathcal C_{\x_0}, \Omega_0)$.  The curves $\mathcal C_\x$ are  of the form \eqref{hyperx} where $u_j(\x)$ are real by Corollary \ref{cor_real}. Therefore, the deformed curves are real and the homology basis of $\mathcal C_{\x_0}$ induces naturally a canonical homolgy basis on $\mathcal C_\x$ with $\x\in\mathcal X\cap \mathbb R^g$ satisfying relations $\rho a_j=a_j$, $\rho b_j=-b_j$ with respect to the antiholomorphic involution. The differentials  $\Omega^{(1)}_{P_\infty}$ on $\mathcal C_\x$ are normalized with respect to this homology basis.  Given that this family is isoperiodic, the Hill curve condition \eqref{eq:hill} holds for $\Omega_0=\Omega^{(1)}_{P_\infty}$ on all the curves $\mathcal C_\x$  for $\x\in\mathcal X\cap \mathbb R^g$.  The  point $\hat P$ remains fixed at infinity, while the branch point $\lambda({\hat Q})$ remains real and in the same real oval as $\hat P$. Thus, the half-integer valued vector $\Delta$ remains fixed.  This  means that the deformed curves generate real solutions \eqref{eq:solSGreal} of the sine-Gordon equation. Since \eqref{eq:periodicity} is satisfied with the same  $T>0$ and $m_j=0$, $j=1, \dots, g,$ for all  the curves $\mathcal C_\x, \; \x\in\mathcal X\cap \mathbb R^g$, as a consequence, solutions $\hat u_\x(X, Y)$ constructed by \eqref{eq:solSGreal} from $\mathcal C_\x$ are  periodic in $X$ with period $T>0$. This proves existence of a continuous family of real periodic in $X$ solutions to the sine-Gordon equation sharing the same period. Equations  \eqref{umxnxk}, \eqref{umxkxk} of Theorem \ref{thm_main} are satisfied since the family $(\mathcal C_\x, \Omega^{(1)}_{P_\infty})$ is isoperiodic, that is since the $b$-period vector $U$ is preserved over the family.
$\Box$

\bigskip

\section{Isoperiodic deformations and comb regions}
\label{sec_CR}

As before, we consider a family $\mathcal C_\x$ of genus $g$ hyperelliptic curves given by equation \eqref{hyperx}. The hyperelliptic coverings of the $\lambda$-sphere have ramification points $P_{\infty}$ over $\lambda=\infty$ along with $P_0,$ $\puj$ and $\pxj$, $j=1, \dots, g$.  We assume $0<u_1<x_1<u_2<\dots < x_{j-1}<u_j<\dots<x_g<\infty$.
Here again $\{\a_1,\dots, \a_g; \b_1, \dots, \b_g\}$ is a fixed canonical homology basis on the surfaces $\mathcal C_\x$ for $\x \in \Sx$, such that   $\a_j$ goes around the points  $\pxj$ and $\puj$ for $j=1,\dots, g$,  while $\b_j$ goes
around the points $P_0,\dots,  P_{u_j}$  for $j=1,\dots, g.$   The domain $\Sx$ is chosen as in Section \ref{sect_iso}. The normalized differential $\Omega_0=\Omega_{P_\infty}^{( 1)}$ is defined by  \eqref{eq:normsecond}.

Define a comb region $\mathcal {CL}$ in  the complex upper half plane $\mathbb H$, as follows:
$$
\mathcal {CL}=\{w\in\mathbb C|\, \Im w>0, \;0\le \Re w\le \infty\}\setminus\bigcup_{j=1}^{g}\{w|\, \Re w = q_j, \; 0\le \Im w\le h_j\}
$$
with some given positive real values $q_1,\dots, q_{g}$ and $h_1, \dots, h_g$.

Differential $\Omega_0=\Omega_{P_\infty}^{(1)}$ having a second order pole at $P_\infty$ on the curve \eqref{hyperx} as its only singularity, can be defined as a polynomial of degree $g$ times the holomorphic differential \eqref{phi}. The coefficients of the polynomial are functions of $\x$ and $u_j(\x)$ ensuring the correct normalization. Denote the zeros of the polynomial by $\xi_j$, $j=1, \dots, g$.   We construct  a Schwarz--Christoffel map $\Theta$, following \cite{MO1975}, as a conformal map $\Theta: \mathbb H\to \mathcal {CL}$:

\begin{equation}
\label{theta}
\Theta(z)=\frac{\i}{2}\int_{0}^z\frac{\prod_j^g(z-\xi_j)\,dz}{\sqrt{z\prod_{j=1}^g(z-x_j)(z-u_j)}},
\end{equation}
mapping   the point at infinity to the point at infinity.

The comb region is a vertical semi-strip with a finite number of  vertical slits where the {\it basis} of the comb region is the semi-infinite ray $[0,\infty)\subset \mathbb R=\{z|\Im z=0\}$ together with the marked points  $q_j$ on it. The  values $q_j$ correspond to the $b$-periods of the differential $\Omega_0$ (up to a constant factor), see \cite{MO1975}. The point  $q_j+\i h_j$ is the $\Theta$-image  of the unique zero $\xi_j$  of the polynomial defining $\Omega_0$ belonging to the interval with the endpoints $\{x_{j}, u_{j}\}$, $j=1, \dots, g$, \cite{MO1975}.

\smallskip

\begin{proposition}\label{prop:CR} Consider  the differential $\Omega_0$  defined by \eqref{Omega} with $\alpha=0$ on a  genus $g$  real hyperelliptic curve $\mathcal C_\x$  with the branch points $\{0, \infty, x_j, u_j|j=1, \dots, g\}$ of the covering $\lambda$ satisfying the ordering conditions stated in this section. Consider  the corresponding  map $\Theta$ \eqref{theta} from the upper half-plane to the comb region defined by positive values $q_j, \;h_j$, $j=1, \dots, g$,   with $q_j$ being the $\Theta$ images of $\{u_j|j=1, \dots, g\}$.  Apply an isoperiodic deformation on the pairs $(\mathcal C_\x, \Omega_0)$ with $\x$ varying in some small neighbourhood $\x\in\mathcal X\subset \mathbb R^g.$ Let $\Theta(\x)$ be the conformal map \eqref{theta} for every $\x\in\mathcal X$. Then, in the image of the $\Theta(\x)$ for $\x\in\mathcal X$, the points $q_j$ remain  independent of $\x$, i.e. the basis of the comb is invariant under the isoperiodic deformations. The deformations
apply only vertically by varying $h_j$, the imaginary parts of the tips of the vertical slits of the comb, that are the imaginary parts of the $\Theta$-images of the zeros of the differential $\Omega_0$, along the vertical rays,  based at $q_j$, for each $j=1,\dots, g$.
Such a family of conformal maps of the comb regions with a finite number of slits and with a semi-infinite basis of the comb, that are obtained by varying $h_j$ vertically are governed by equations \eqref{umxnxk}, \eqref{umxkxk} as well as \eqref{derivative}, as equations for the real-valued functions $u_1(\x), \dots, u_g(\x)$ from \eqref{hyperx}.
\end{proposition}

{\it Proof.} The values $q_j$ correspond to the $b$-periods of  the differential $\Omega_0$, by construction and properties of the Schwarz--Christoffel map $\Theta$, as explained in \cite{MO1975}. The isoperiodic deformations  $(\mathcal C_{\x}, \Omega_0)$ preserve the periods of the deformed differential, thus they preserve the basis of the comb. Therefore, during the deformation, the tips of the vertical slits, defined by $h_j$ deform vertically, along the slit. The isoperiodic deformations are governed by equations \eqref{umxnxk}, \eqref{umxkxk} for the real-valued functions $u_1(\x), \dots, u_g(\x)$ from \eqref{hyperx}, according to Theorem \ref{thm_main} and Corollary \ref{cor_real}, and satisfy \eqref{derivative}. $\Box$

\

Thus, Proposition \ref{prop:CR} provides an explicit rectification of isoperiodic deformations, seen as vertical deformations of the teeth of the comb, with an unchanged base of the comb,  together with the exact equations which govern deformations of the $\Theta$-preimages of the base of the obtained comb regions.

\begin{remark}\label{rem:Loewner}
In the case $g=1$, the comb regions have one vertical slit. The vertical deformations of that slit are governed by differential equation \eqref{ode} from Theorem \ref{thm_ode}. Thus, equation \eqref{ode}  is related to a version of the Loewner differential equation, which governs deformations of conformal maps induced by deformations of regions with one slit.   Therefore, we can  see equations \eqref{umxnxk}, \eqref{umxkxk} of Theorem \ref{thm_main} as a sort of PDE generalizations related to multi-slit regions. For more about Loewner differential equation, see \cite{Loewner, Kufarev, Law, Duren, GN}.
\end{remark}

\bigskip

{\bf Acknowledgements.}  V.D. acknowledges with gratitude the Simons Foundation grant no. 854861 and the Science Fund of Serbia grant IntegraRS. V.S. gratefully acknowledges
support from the Natural Sciences and Engineering Research Council of Canada through a Discovery grant and from the University of Sherbrooke.

\end{document}